\documentclass[a4paper]{article}

\usepackage[plainpages=false, colorlinks=true, 
            linkcolor=black, urlcolor=black, citecolor=black]{hyperref}
\usepackage{geometry}

\usepackage{tikz}
\usepackage{overpic}
\usepackage{amsmath, amssymb, amsthm,  amsfonts, mathrsfs}
\usepackage{multirow}
\usepackage{hhline}
\usepackage{array}
\usepackage{hyperref}
\hypersetup{colorlinks,linkcolor={blue},citecolor={red!70!black},urlcolor={red}}  
\usepackage{color}
\usepackage{booktabs}
\usepackage{tikz}
\usepackage{bm}
\usepackage{scrextend}
\usepackage{arydshln}
\usepackage {url}

\usepackage[font=scriptsize, labelfont=bf]{caption}
\usepackage{graphicx} 
\usepackage{epstopdf} 
\usepackage{algpseudocode}
\usepackage{array,  tabularx}
\usepackage{algorithm, algcompatible}
\usepackage{booktabs} 
\usepackage{paralist} 
\usepackage{verbatim} 
\usepackage{subfig} 
\usepackage{xcolor, marginnote, enumitem}
\usepackage[numbers]{natbib}
\numberwithin{equation}{section}
\theoremstyle{theorem}
\newtheorem{lemma}{Lemma}
\newtheorem{theorem}{Theorem}
\newtheorem{proposition}{Proposition}

\theoremstyle{remark}
\newtheorem{remark}{Remark}

\theoremstyle{definition}
\newtheorem{definition}{Definition}

\DeclareMathOperator{\gph}{gph}
\DeclareMathOperator{\dist}{dist}

\DeclareMathOperator{\dom}{dom}

\DeclareMathOperator*{\argmin}{arg\, min}

\newcommand{\grad}{\nabla}
\newcommand{\bb}{\mathbb}

\newcommand{\lookUp}[1]{}

\algnewcommand\INPUT{\item[\textbf{Input:}]}%
\algnewcommand\OUTPUT{\item[\textbf{Output:}]}%

\DeclareMathAlphabet{\mathbfit}{OML}{cmm}{b}{it}

\begin{document}

\title{Preconditioned Algorithm for Difference of Convex Functions with applications to Graph Ginzburg-Landau Model}

\author{
Xinhua Shen\thanks{Institute for Mathematical Sciences, 
		Renmin University of China,  China.
		Email: \href{mailto:shenxinhua@ruc.edu.cn}{shenxinhua@ruc.edu.cn}.} 
			\quad 
Hongpeng Sun\thanks{Institute for Mathematical Sciences,
Renmin University of China, China.
Email: \href{mailto:hpsun@amss.ac.cn}{hpsun@amss.ac.cn}.} \quad 
 Xuecheng Tai\thanks{NORCE Norwegian Research Centre, Nygardstangen, NO-5838 Bergen, Norway. Email: \href{xtai@norceresearch.no}{xtai@norceresearch.no}.} 
}

\maketitle
\begin{abstract}
In this work, we propose and study a preconditioned framework with a graphic Ginzburg– Landau functional
 for image segmentation and data clustering by parallel computing. Solving nonlocal models is usually challenging due to the huge computation burden. For the nonconvex and nonlocal variational functional, we propose several damped Jacobi and generalized Richardson preconditioners for the large-scale linear systems within a difference of convex functions algorithms framework.  They are efficient for parallel computing with GPU and can leverage the computational cost. Our framework also provides flexible step sizes with a global convergence guarantee. Numerical experiments show the proposed algorithms are very competitive compared to the singular value decomposition based spectral method. 
\end{abstract}

\paragraph{Key words.}{}
graph model, Ginzburg-Landau functional, difference of convex functions algorithms, damped Jacobi preconditioner, Richardson preconditioner, parallel computing,  image segmentation, data clustering



\section{Introduction}
\label{sec:intro}
Graph models provide powerful frameworks and tools with unified representation for image processing and analysis \cite{LG}. Especially, the nonlocal graph models are widely used for data clustering \cite{BF,EBAT,ZS} and computer vision tasks including optical flow estimates \cite{KKr,RBP}, image denoising \cite{GO2}, and image segmentation \cite{SM,SB,Tang2019}. Detailed analysis and mathematical interpretations can be found in  \cite{BCM,GO1,GO2}. Besides, nonlocal graph models are also developed as regularization techniques for inverse problems \cite{BCH,PBC} including inverse acoustic scattering \cite{CK}. For its recent application with machine learning, we refer to \cite{CPKMY}. It is believed that the nonlocal graph model can capture certain important global information of images or data compared to local models \cite{SM}.

Inspired by the recent development of Ginzburg-Landau (GL) functional-based nonlocal variational framework \cite{BF, GMBFP, LB, MBYL}, we propose a preconditioned difference of convex functions algorithm (DCA) for image segmentation and data clustering. Our preconditioners aim to solve the extremely large-scale linear system produced by the nonlocal graph model, which is usually very difficult and challenging to compute. The specially designed damped Jacobi and generalized Richardson preconditioners are very efficient. These preconditioners have inherited advantages with parallelization. Global convergence of the proposed algorithm is guaranteed. Numerical tests show that preconditioned iterations can provide high-quality image segmentation and data clustering with as low as four iterations. Our algorithms offer a good alternative to circumvent the difficulties related to the high computational cost.  We also prove the global convergence of the proposed preconditioned DCA with the novel extensions of the Kurdyka-\L ojasiewicz (KL) analysis for nonconvex optimization.

Compared to other existing techniques based on singular value decomposition (SVD) and spectral method \cite{BF,GMBFP,LB,ME,MTB}, our framework has the following advantages. First,  its structure is especially well-suited for parallel computing. Numerical experiments show that the proposed preconditioned framework can be more than  10 times faster compared with the time for image segmentation tasks with large images. Although the Nystr\"{o}m method \cite{FBCM}  only needs to compute a few eigenvalues and eigenvectors, the number of needed eigenvalues and eigenvectors are still very high even with a very low sampling rate for large images or big data sets. Furthermore, the segmentation or data clustering quality heavily depends on the sampling rate. For our preconditioned framework, its efficiency can be greatly improved with preconditioned Jacobi or Richardson iteration through a sparse window via GPU instead of the computation of SVD.  Second, we proposed two different techniques to compute the weights for the graph Laplacian. Tests show that they can handle different situations, like the segmentation task for images with disconnected components. Third, compared to some existing convex splitting frameworks with conditional stability for step sizes as in \cite{BF,LB,MTB}, 
we introduced the preconditioned DCA framework that is stable for any positive step size.  
While the step size goes into infinity, we can get an efficient preconditioned DCA without step size.

The remaining of this work is organized as follows. Section \ref{sec:intro:graphlap:GL}  gives an introduction to the graph Laplacian-based Ginzburg-Landau functional. In section \ref{sec:pre:dca:main}, we build the proposed preconditioned DCA framework for arbitrary positive step sizes. With a detailed analysis of the graph Laplacian, two kinds of feasible Jacobi preconditioners and one feasible generalized Richardson preconditioner are introduced for the unnormalized and normalized graph Laplacian cases.  In section \ref{sec:global:conv}, we prove the global convergence of the proposed algorithms using the Kurdyka-\L ojasiewicz technique. A discussion of the local convergence rate is also supplied. In section \ref{sec:num}, detailed numerical experiments for image segmentation and data clustering are presented. We also introduce two novel nonlocal graph Laplacian matrices for image segmentation. Numerical examples show the efficiency of the proposed framework. Finally, some concluding remarks are given in section \ref{sec:conclusion}.

\section{The graph Laplacian based Ginzburg-Landau functional}\label{sec:intro:graphlap:GL}
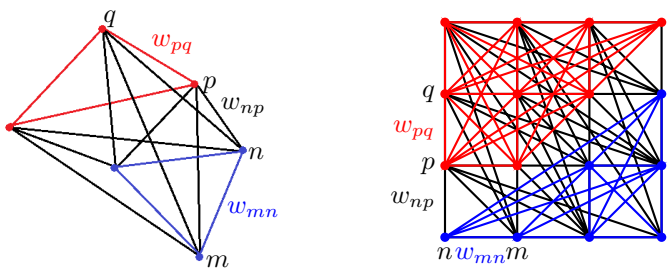
\begin{figure}[]
	\centering 
	\begin{minipage}{0.45\linewidth}
		\qquad \qquad
		\begin{overpic}[width = 0.6\textwidth]{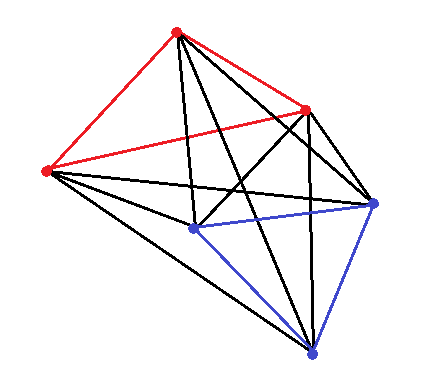}
			\put(76,5){\color{black}{$m$}}
			\put(90,41){\color{black}{$n$}}
			\put(75,65){\color{black}{$p$}}
			\put(42,87){\color{black}{$q$}}
			\put(83,23){\color{blue}{$w_{mn}$}}
			\put(81,58){\color{black}{$w_{np}$}}
			\put(58,79){\color{red}{$w_{pq}$}}
		\end{overpic}
	\end{minipage} 
	\begin{minipage}{0.45\linewidth}
		\begin{tikzpicture}[thick,scale=0.95]
		\draw[step=1cm] (0,0) grid (3,3);
		\draw[black] (0,1) -- (1,0) (1,1) -- (2,0) (1,2) -- (2,1) (2,2) -- (3,1) (2,3) -- (3,2)
		(0,2) -- (1,0) (0,3) -- (1,0) (0,1) -- (2,0) (0,3) -- (2,0) (1,3) -- (2,0) (1,2) -- (2,0) (0,2) -- (2,1) (1,3) -- (2,1) (0,1) -- (3,0) (1,3) -- (3,0) (2,3) -- (3,0) (0,2) -- (3,1) (0,3) -- (3,1) (1,2) -- (3,1) (2,3) -- (3,1) (0,3) -- (3,2) (1,3) -- (3,2);
		\draw[blue] (0,0) -- (3,0) -- (3,2) (2,0) -- (2,1) -- (3,1);
		\draw[blue] (2,0) -- (3,1) -- (1,0) -- (3,2) -- (0,0)
		(0,0) -- (2,1) (0,0) -- (3,1) (2,0) -- (3,2)
		(3,0) -- (2,1);
		\draw[red] (3,3) -- (0,3) -- (0,1) -- (1,1) -- (1,3) (0,2) -- (2,2) -- (2,3);
		\draw[red] (0,2) -- (1,3) -- (0,1) -- (2,3) -- (1,1) -- (3,3)
		(3,3) -- (0,2) (3,3) -- (1,2) (3,3) -- (0,1)
		(2,3) -- (0,2) (2,2) -- (0,1)
		(0,2) -- (1,1) -- (0,3) -- (1,2) 
		(0,3) -- (2,2) (1,3) -- (2,2) ;

		\draw[fill,blue] (0,0) circle [radius=0.05]
		node[below, black] {$n$};
		\draw[fill,red] (0,1) circle [radius=0.05]
		node[left, black] {$p$};
		\draw[fill,red] (0,2) circle [radius=0.05]
		node[left, black] {$q$};
		\draw[fill,blue] (1,0) circle [radius=0.05]
		node[below, black] {$m$};
		\draw[fill,red] (0,3) circle [radius=0.05]
		node[below, red] {};
		\draw[fill,black] (0,0.5) circle [radius=0]
		node[left, black] {$w_{np}$};
		\draw[fill,black] (0,1.5) circle [radius=0]
		node[left, red] {$w_{pq}$};
		\draw[fill,black] (0.5,0) circle [radius=0]
		node[below, blue] {$w_{mn}$};
		
		\draw[fill,red] (1,1) circle [radius=0.05]
		node[below, red] {};
		\draw[fill,red] (1,2) circle [radius=0.05]
		node[below, red] {};
		\draw[fill,red] (1,3) circle [radius=0.05]
		node[below, red] {};
		\draw[fill,red] (2,2) circle [radius=0.05]
		node[below, red] {};
		\draw[fill,red] (2,3) circle [radius=0.05]
		node[below, red] {};
		\draw[fill,red] (3,3) circle [radius=0.05]
		node[below, red] {};
		\draw[fill,blue] (3,2) circle [radius=0.05]
		node[below, red] {};
		\draw[fill,blue] (3,1) circle [radius=0.05]
		node[below, red] {};
		\draw[fill,blue] (3,0) circle [radius=0.05]
		node[below, red] {};
		\draw[fill,blue] (2,1) circle [radius=0.05]
		node[below, red] {};
		\draw[fill,blue] (2,0) circle [radius=0.05]
		node[below, red] {};
		\end{tikzpicture}
	\end{minipage}
	\caption{The general undirected graph (left) and regular undirected graph (right)}
	\label{tab:my_label:undirected:graph}
\end{figure}

Let us first turn to the graph structure, which is important for the construction of the graph Laplacian operator for these models. It contains the relationships among the image's pixels or data. We can regard the image or data as an undirected weighted graph $G = (V,  E)$, where $V $ is the set of nodes corresponding to the image's pixels or the data points,  and $E \subset V \times V$ is the set of edges \cite{LG,ME}. The edge weighting is a function $w : E \rightarrow \mathbb{R}$  corresponding to the weights set $\{w_{ij}\}$ with $w_{i,j} = w(e_{i,j})\geq 0$ and the edge $e_{i,j} \in E$ \cite{LG}; see Figure \ref{tab:my_label:undirected:graph}. 

\subsection{Ginzburg-Landau functional}
Let us directly focus on the following discrete graph Ginzburg-Landau (GL) functional for image segmentation and data clustering
\begin{equation}\label{eq:gl:functional:graph}
GL(u) := \frac{\epsilon}{2}\langle u, Lu\rangle + \frac{1}{\epsilon}\sum_i W(u(i)).
\end{equation}
Here $W(u) = \frac{1}{4}(u^2 - 1)^2$ is the double-well potential,  $\epsilon$ is the small diffuse parameter, and the graph Laplacian $L$ with symmetric weights $w_{i,j}=w_{j,i}\geq 0$ is defined as follows
\begin{equation}
\langle u, Lu\rangle := \frac12\sum_{ij}w_{ij}(u(i)-u(j))^2.
\end{equation}
Here and later, $u(i)$ is the labeling variable for the $i$-th pixel or data point.  It has been proven that the minimizer of the GL functional \eqref{eq:gl:functional:graph}  converges to a binary label $u \in \{ 1,-1 \}$. While the small diffuse parameter $\epsilon$ tends to 0, the GL functional will gamma-converge to the total variation (TV) seminorm, which is also a graph cut function that is powerful for segmentation or clustering  \cite{BF}.     
We refer to \cite{BF} for more discussions on the corresponding mechanism. 
The graph Ginzburg-Landau functional \eqref{eq:gl:functional:graph} originated from the following continuous Ginzburg-Landau (GL) functional \cite{BF, LB},
\begin{equation}\label{eq:gl:functional}
GL_{c}(u) := \frac{\epsilon}{2}\int|\grad u|^2d\sigma + \frac{1}{\epsilon}\int W(u)  d\sigma.
\end{equation}
For certain  domains $\Omega \subset \mathbb{R}^n$ and small parameter $\epsilon$, there exist nonconstant local minimizers $u^{\epsilon}$ satisfying $u^{\epsilon} \approx \pm 1$ except in a thin transition region \cite{kohn_sternberg_1989}. The double
well potential $W$ in \eqref{eq:gl:functional}  will force $u$ to go to either one or minus one. The $H^1$ term forces $u$ to
have some smoothness and remove sharp jumps between the two minima of $W$.  This leads to a thin transition region with $\epsilon$ scale \cite{BF}. The GL functional is widely employed for multiscale phase-field simulations\cite{Ey}.

By adding a weighted quadratic fidelity term to \eqref{eq:gl:functional:graph}, the following minimization problem can be used for many applications including image segmentation or data clustering:
\begin{equation}\label{eq:functional:GL}
u=\argmin_{u} F(u), \quad    F(u): = GL(u) + \frac{\eta}{2} ||u-y||_{\Lambda}^2.
\end{equation}
Here, $||u-y||_{\Lambda}^2 :=\langle u-y, \Lambda (u-y)\rangle = \sum_{i} \Lambda_{i}(u(i)-y(i))^2$ and  $\Lambda \geq 0$ is a diagonal matrix and the values $y(i)$ are assumed to be known. In applications, $\Lambda_{i} = 1$ if $i \in Z'$ and 0 otherwise with $Z'$ being the set of pixels or data points with known prior labels. The balancing parameter $\eta$ is positive, i.e., $\eta>0$. Let's  denote
\[
\bb{W}(u) = \frac{1}{\epsilon}\sum_i W(u(i)).
\]
We thus can write $F(u)$ as follows
\begin{equation}\label{eq:discrete:fu}
F(u) = \sum_{ij}\frac{\epsilon}{2}w_{ij}(u(i)-u(j))^2 +\bb{W}(u) + \frac{\eta}{2}\sum_i\Lambda(i)(u(i)-y(i))^2.
\end{equation}
The weight $w_{i,j}$ is critical for the graph GL model. Now let us turn to the construction of the graph Laplacian and the weights.
\subsection{The Graph Laplacian and the Weight}

We compute the weight $w_{ij}$ between pixel $i$ and pixel $j$ (or the data point $i$ and data point $j$) by combining a feature similarity term $K(i,j)$ and a proximity term $N(i,j)$ as follows
\begin{equation*}
w_{ij} =K(i,j) \cdot N(i,j).
\end{equation*}
Here, $K(i,j)$ is a Gaussian kernel function, and $N(i,j)$ is an indicator function.

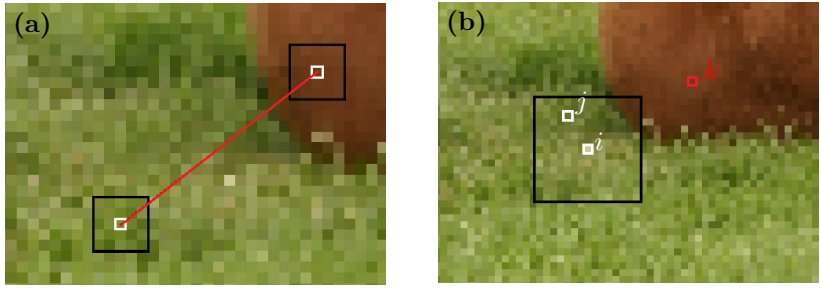
\begin{figure*} [t!]
	\centering
	\begin{minipage}{0.34\linewidth}
		\vspace{0.47cm}
		\begin{overpic}[width = 1\textwidth]{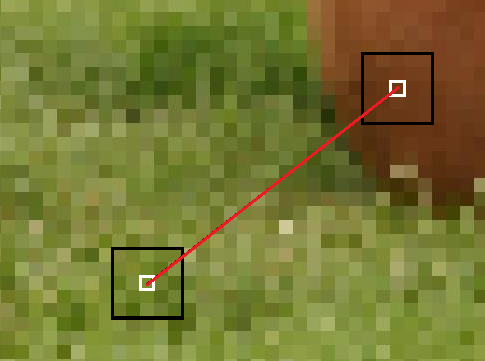}
			\put(2,67){\color{black}{\textbf{(a)}}}
		\end{overpic}
		\label{fig:patch_similarity}
	\end{minipage}\qquad
	\begin{minipage}{0.34\linewidth}
		\begin{overpic}[width = 1\textwidth]{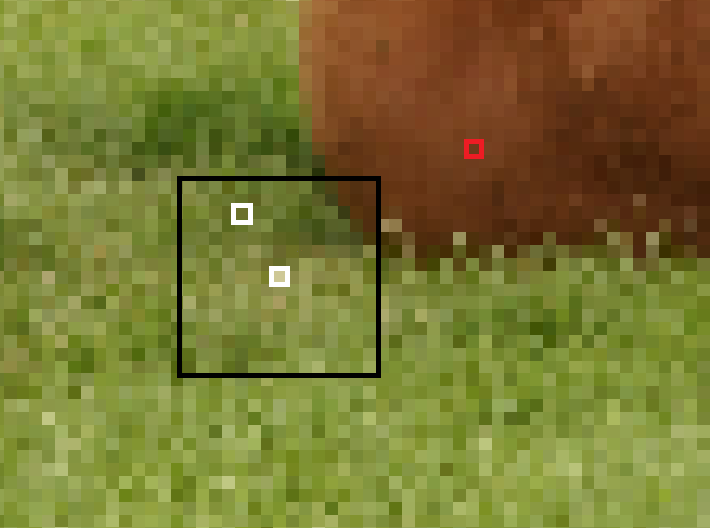}
			\put(2,67){\color{black}{\textbf{(b)}}}
			\put(41,35){\color{white}{$i$}}
			\put(36,46){\color{white}{$j$}}
			\put(70,54){\color{red}{$k$}}
		\end{overpic}
	\end{minipage}
	\caption{Figure {\bf{a}} shows how to calculate $K(i,j)$ in \eqref{eq:kij} between two patches. Figure {\bf{b}} shows the search window and the proximity in a search window. While calculating the weight inside the window, $w_{i,j} \neq 0$  and $w_{i,k}=0$ since the pixel $k$ is not in the search window. }
	\label{patch_and_window} 
\end{figure*}

Now let's first show the ways to compute $K(i,j)$. For images, we define two different neighbors of a pixel including {\it patch} and {\it window}, to introduce the computing of weights in the graph. {\it Patch} is a $5 \times 5$ neighbor (Fig.\ref{patch_and_window}a) to estimate the similarity between two pixels.  {\it Window} is a neighbor (Fig.\ref{patch_and_window}b) to determine the proximity between two pixels. In general, the size of a {\it window} is $15 \times 15$ pixels otherwise specified in section \ref{sec:num}.
Define the color information of pixel $i$ is $C_i = (C^R_i, C^G_i, C^B_i)$, and the Euclidean distance $||C_i - C_j||_2$ shows the difference between two pixels in color. To describe the similarity between two pixels exactly, we consider their patches. Computing a weighted sum and mapping it to $(0,1]$ by the Gaussian kernel function, we can get the formula of similarity term $K(i,j)$,
\begin{equation}\label{eq:kij}
K(i,j) = \exp\left(-\frac{\sum_{n=1}^{n_0}\alpha_n||D_{in} - D_{jn}||_2^2}{\sigma^2}\right).
\end{equation}
For images, $D_{in} = (C_{in}^R, C_{in}^G, C_{in}^B)$ represents the color information about the $n$-th neighbor of pixel $i$ and $\alpha_n$ depends on the distance between the $n$-th neighbor and the center pixel in the patch. As in Figure \ref{patch_and_window}a, $(D_{i0}, \ldots, D_{in_0})$ and $(D_{j0}, \ldots, D_{jn_0})$ are within the same order in the corresponding $i$-th and $j$-th patch. $n_0$ denotes the patch size (e.g., $n_0=25$ for $5 \times 5$ neighbor in Fig.\ref{patch_and_window}a).   For clustering, $D_{i} = (D_{i1}, \ldots,D_{in_0})$ represents the position of $i$-th point in the feature space, $\alpha_n$ sets to 1 and $n_0$ represents the dimension of the feature space. Here $\|D_i-D_j\|_2^2=\sum_{n=1}^{n_0}||D_{in} - D_{jn}||_2^2$ and $\|D_i-D_j\|_2$ is actually the distance between the $i$-th data point and the $j$-th data point in the feature space.



Now, let us turn to compute the proximity term $N(i,j)$.  For images, 
the proximity term determines whether a pixel relates to another pixel by the window. As in Figure \ref{patch_and_window}b, the box with a black boundary is the window of pixel $i$. It contains pixel $j$, but not
$k$. Assuming the locations of pixel $i$, pixel $j$ are $P_i = (x_i,y_i)$, $P_j = (x_j,y_j)$, we define the proximity terms as
\begin{equation*}
N_s(i,j) = \begin{cases}
1,      \quad   \max\{|x_i-x_j|,|y_i-y_j|\} < \text{dist},\\
0,      \quad \text{others},
\end{cases} \ \  
N_c(i,j) = \begin{cases}
1, \quad &\|D_i - D_j\|_2 \leq \|D_i - D_{n_l}\|_2,\\
0, \quad &\text{other}.
\end{cases}
\end{equation*}
Here $N(i,j)=N_s(i,j)$ for images with `$\text{dist}$' determined by the window. For clustering, $N(i,j)=N_c(i,j)$ with predetermined $n_l$. Moreover, we employ the KNN for searching nearest neighbors, i.e., the $k$-nearest neighbors algorithm based on the Euclidean distance $\|\cdot\|_2$. We refer to \cite{MPL} for the implementations of KNN with the source code. 

Now, we obtain a weight matrix $W=(w_{i,j})$ of the graph $G$ and introduce the following unnormalized graph Laplacian,
\begin{equation*}
L= D - W.
\end{equation*}
We will also introduce the following normalized graph Laplacian
\begin{equation*}
L_s=I-D^{-\frac{1}{2}}WD^{-\frac{1}{2}}
\end{equation*}
where $D$ is a diagonal matrix $d_{ii} = \sum_j w_{ij} >0$.
$L$ and $L_s$ are based on an assumption: each node is connected to any other node in the window and the corresponding weight is not 0. It can be seen that the sparsity of the graph Laplacian $L$ or the normalized graph Laplacian $L_s$ depends on the weight matrix $W$, especially the proximity term $N(i,j)$. The increase in sparsity will increase the computation cost.

For the parameter $\sigma$ in the Gaussian kernel $K(i,j)$ as in \eqref{eq:kij}, we choose fixed bandwidth $\sigma$ to be on the order $\log(N_d)+1$ for image segmentation, where $N_d$ is the number of vertices (see section 2.2.1 of \cite{MTB}). For data clustering, we employ the local scaling weight
with adaptive bandwidth introduced by Zelnik-Manor and Perona \cite{MLP} (see formula (2.9) of \cite{BF} or (2.12) in \cite{MTB}).

We now turn to the proposed preconditioned DCA for the graph Ginzburg-Landau model. 
\section{Preconditioned difference of convex functions algorithm} \label{sec:pre:dca:main}
Let's start with minimizing the nonlinear functional \eqref{eq:functional:GL}. First, write $F(u)$ as the difference of convex functions
\begin{equation}\label{eq:dca:dc}
F(u) = P_1(u)-P_2(u),  
\end{equation}
with 
\begin{equation}\label{eq:dc:functional}
P_1(u) :=\frac{\epsilon}{2}\langle u, Lu \rangle +  \frac12 \eta ||u-y||_{\Lambda}^2 + \frac{c}{2}\|u\|_2^2, \quad P_2(u) := \frac{c}{2}\|u\|_2^2 - \bb{W}(u).
\end{equation}
The format of \eqref{eq:dca:dc} decomposes a complex nonlinear functional into one simple convex term and the other term. Different from the previous work \cite{GMBFP}, we put the fidelity term $\frac{1}{2}\eta\|u-y\|_{\Lambda}$ on the $P_1(u)$ with the Laplacian term. This change can bring convenience for computing and testing.
Here, the constant $c>0$ is chosen such that $P_2$ is convex on $u$. 
The standard DCA iterations then follow, for $t=0, 1, \ldots$,
\begin{equation}\label{eq:orig:dca}
u^{t+1}:=\argmin_{u} P_1(u) - \langle\xi^{t}, u\rangle,  \quad \xi^t \in  \partial P_2(u)|_{u=u^t}.
\end{equation}
It turns the nonlinear and nonconvex minimization problem \eqref{eq:functional:GL} to a sequence of convex minimization problems \eqref{eq:orig:dca} by the linearization of  $P_2$. With $\xi^t = cu^t - \grad \bb W(u)|_{u=u^t}$, we need to solve the following linear equation during each DCA iteration as in \eqref{eq:orig:dca}
\begin{equation}\label{eq:dca:linear:system:ori}
Tu^{t+1}=b^t, \quad T:=    \epsilon L + \eta \Lambda +cI, \ \ b^t :=  \eta \Lambda y +\xi^t.
\end{equation}
It is very challenging to solve such kind of large linear systems.  In \cite{DS}, a preconditioned DCA framework is proposed for dealing with this kind of problem. Global convergence can be guaranteed with any finite feasible preconditioned iterations for the linear system appeared in \eqref{eq:orig:dca}. The preconditioned DCA can be obtained by following weighted and proximal DCA
\begin{equation}\label{eq:pdca:framework}
u^{t+1}:=\argmin_{u} P_1(u) - \langle\xi^{t} , u\rangle+ \frac{1}{2}\|u-u^k\|_{M_0}^2,  \quad \xi^t  \in \partial P_2(u)|_{u=u^t}, \tag{pDCA}
\end{equation}
where $M_0$ is a positive definite and self-adjoint operator (or matrix). The  proximal term $\frac{1}{2}\|u-u^k\|_{M_0}^2$ can help introduce preconditioned iterations \cite{DS}. The weighted norm $\|\cdot\|_{M_0}$ is defined as follows
\begin{equation}
||x||_{M_0}^2 :=\langle x , M_0x\rangle, \quad \text{with} \ \   \langle x, y \rangle_{M_0} := \langle x, M_0y \rangle.
\end{equation}
However, the step size is not considered in \cite{DS}. It is discussed in \cite{Ey} in the framework of the discretization of the gradient flow of minimal functional. Let's turn to $F(u)$ in \eqref{eq:functional:GL} with a difference of convex functional \eqref{eq:dc:functional}. Inspired by the semi-implicit splitting method as in \cite{LB} and  the linearly stabilized splitting scheme for the Cahn-Hilliard equation as in \cite{Ey},
the DCA for the graph Ginzburg-Landau model with step size $k$ can be written as follows 
\begin{equation}\label{eq:dca:stepsize}
\frac{u^{t+1}-u^t}{k} =  - (\nabla P_1(u)|_{u=u^{t+1}}- \xi^{t}),   \quad \xi^t  \in \partial P_2(u)|_{u=u^t}.
\end{equation}
With direct calculation, the linear system for $u^{t+1}$ in \eqref{eq:dca:stepsize} becomes
\begin{equation}\label{eq:orig:equation:with:k}
T_k u^{t+1} = b_k^t, \quad T_k = I + kT, \quad b_k^t = u^t + kb^t,
\end{equation}
where $T$ and $b^t$ are the same as in \eqref{eq:dca:linear:system:ori}.

It is convenient to employ varying step sizes. As shown in \cite{LB}, a monotone decrease of the energy $F(u)$ can be obtained under certain step size constraints.
We will introduce step sizes and preconditioners simultaneously within the framework of preconditioned DCA through the following lemma.

\begin{lemma}\label{lemma:preconditioned_DCA:stepsize}
	With appropriately positive constant $\kappa$, linear and bounded operator $M$ and setting $M_0 = \frac{1}{k}I + \frac{\kappa}{k}M$ such that $M_0 \geq L_0I$ with positive constant $L_0$,  the iteration \eqref{eq:pdca:framework} can be reformulated as the following classical preconditioned iteration for \eqref{eq:orig:equation:with:k}
	\begin{equation}\label{eq:orig:equation:with:k:pre}
	u^{t+1}:=u^t + \mathbb{M}_{k,\kappa}^{-1}(b_k^t - T_k u^t).
	\end{equation}
	Here,
	\begin{equation}\label{eq:btk:utk}
	b_k^t = u^t + k(cu^t - \grad \bb W(u)|_{u=u^t} +\eta\Lambda y), \ \  T_k =  I + k(\epsilon L + \eta\Lambda + cI), \ \ \mathbb{M}_{k,\kappa} = \kappa M + T_k.
	\end{equation}
\end{lemma}
\begin{proof}
	For the preconditioned DCA, with the same notations in \eqref{eq:pdca:framework}, we have
	\begin{equation}\label{eq:pdca:framework:stepsize}
	u^{t+1}:=\argmin_{u} P_1(u) - \langle\xi^{t} , u\rangle+ \frac{1}{2}\|u-u^k\|_{\frac{I}{k}}^2 + \frac{1}{2}\|u-u^k\|_{\frac{\kappa}{k}M}^2,  \quad \xi^t  \in \partial P_2(u)|_{u=u^t}.   \tag{pDCA-k}
	\end{equation}
	By direct calculation, we arrive at
	\begin{align}
	&\epsilon L u^{t+1} + \eta \Lambda (u^{t+i} - y) + cu^{t+1} + M(u^{t+1}-u^t)/k + (u^{t+1}-u^t)/k= \xi^t  \notag \\
	&\Longleftrightarrow \kappa M(u^{t+1}-u^t) + (I + k(\epsilon L + \eta\Lambda + cI))u^{t+1} = k (\xi^t + \eta \Lambda y) + u^k.
	\end{align}    
	With the notations $b_k^t$ and $T_k$  in \eqref{eq:btk:utk}, we further have
	\begin{equation}\label{eq:pre:step:prox:M}
	\kappa M(u^{t+1} - u^t) + T_k u^{t+1} = b_k^t.
	\end{equation}
	It can be reformulated as the classical stationary preconditioned iteration as follows,
	\begin{align*}
	& (M + T_k)u^{t+1} = (M + T_k)u^t + (b_k^t - T_k^t ) \notag \\
	&\Longleftrightarrow  u^{t+1} = u^t  + (M + T_k)^{-1}(b_k^t - T_k^t ) = u^t + \mathbb{M}_{k,\kappa}^{-1}(b_k^t - T_k^t ), \quad \mathbb{M}_{k,\kappa}^{-1}: = \kappa M + T_k.
	\end{align*}
\end{proof}

Similarly, for the preconditioned DCA without step sizes \eqref{eq:pdca:framework}, we have the following proposition, whose proof is quite similar to Lemma \ref{lemma:preconditioned_DCA:stepsize} and is omitted.
\begin{proposition}\label{prop:nopre}
	With appropriately chosen linear operators $M_0 \geq \delta_0I$ and  constant $\delta_0>0$,  the iteration \eqref{eq:pdca:framework} can be reformulated as the following classical preconditioned iteration for \eqref{eq:dca:linear:system:ori}
	\begin{equation}\label{eq:dca:linear:system:ori:pre}
	u^{t+1}:=u^t + \mathbb{M}^{-1}(b^t - T u^t).
	\end{equation}
	Here,
	\begin{equation}\label{eq:bt:ut}
	b^t = cu^t - \grad \bb W(u)|_{u=u^t} +\eta\Lambda y, \ \  T =  \epsilon L + \eta\Lambda + cI, \ \ \mathbb{M} = M_0 + T.
	\end{equation}
\end{proposition}
Letting $k\rightarrow +\infty$, we have the following connections between the finite step size cases and the infinite step size cases.
\begin{proposition}\label{prop:k:infinity} 
	Assuming the iteration sequence $\{u^t\}$ of proximal DCA \eqref{eq:pdca:framework} is bounded, then the update of $u^{t+1}$ in \eqref{eq:dca:stepsize}  without preconditioners can be obtained by setting $M_0={I}/{k}$ in \eqref{eq:pdca:framework}. Furthermore,  the linear update \eqref{eq:orig:equation:with:k} will converge to \eqref{eq:dca:linear:system:ori}  and the linear update  \eqref{eq:orig:equation:with:k:pre} will converge to   \eqref{eq:dca:linear:system:ori:pre} for $\kappa=k$ when $k\rightarrow +\infty$. 
\end{proposition}
\begin{proof}
	It can be readily checked that while setting $M_0={I}/{k}$ in \eqref{eq:pdca:framework}, we can get the linear update \eqref{eq:orig:equation:with:k} for DCA with step size $k$. 
	For the linear update \eqref{eq:orig:equation:with:k}, it can be reformulated as 
	\[
	(T_k/k)u^{t+1} = b_{k}^t/k \Leftrightarrow (I/k +T)u^{t+1} = u^t/k + b^t \Longrightarrow Tu^{t+1}=b^t
	\]
	while $k\rightarrow +\infty$ with the boundedness assumption of $\{u^t\}$. Similarly, for the linear update \eqref{eq:orig:equation:with:k:pre}, dividing  both sides of the equivalent linear equation \eqref{eq:pre:step:prox:M} by $k$, we obtain
	\[
	\frac{\kappa}{k}M(u^{t+1}-u^t) + T_k/k=b_k^t/k \Longrightarrow  M(u^{t+1}-u^t) + T u^t =b^t.
	\]
	The last equation is obtained by letting $k=\kappa\rightarrow +\infty$ and the boundedness assumption of $\{u^t\}$. It is the same as  \eqref{eq:dca:linear:system:ori:pre} when choosing $M_0=M$ finally.
\end{proof}
The boundedness of $\{u^t\}$ of proximal DCA \eqref{eq:pdca:framework} in Proposition \ref{prop:nopre} will be proved later. 
\begin{remark}
	Inspired by \cite{BF,GMBFP,LB}, instead of the difference of convex functional  in  \eqref{eq:dca:dc}, we can also split $F(u) = P_1(u)-P_2(u)$ as follows, 
	\begin{equation}\label{eq:dc:functional:svd}
	P_1(u) :=\frac{\epsilon}{2}\langle u, Lu \rangle + \frac{c}{2}\|u\|_2^2, \quad P_2(u) := \frac{c}{2}\|u\|_2^2 - \bb{W}(u) -  \frac12 \eta ||u-y||_{\Lambda}^2.
	\end{equation}
	Here, the constant $c>0$ is chosen such that $P_2$ is a convex functional on $u$. The SVD-Nystr\"{o}m based methods \cite{BF,GMBFP}  prefer the difference of convex functional in \eqref{eq:dc:functional:svd} with the constant coefficient linear operator $\epsilon L + cI$ to $T = \epsilon L + cI + \eta \Lambda $ in \eqref{eq:dca:linear:system:ori}. The reason is that the linear operator   $\epsilon L + cI$ is independent of the prior $\eta \Lambda$. However, the proposed various preconditioners in our framework can do preconditioning for $T$ directly with varying coefficient $\Lambda$ efficiently. 
	\label{rem:dca:constantcoefficient}
\end{remark}
Due to the Proposition \ref{prop:k:infinity}, henceforth, we call the DCA with \eqref{eq:dca:linear:system:ori} and  \eqref{eq:pdca:framework} with \eqref{eq:dca:linear:system:ori:pre} as the case of step size infinity compared to the finite step sizes cases.
Now, let us turn to the detailed preconditioners for the preconditioned DCA. 
\subsection{Feasible Jacobi and Richardson Preconditioners}\label{subsec:preconditioner}
In this paper, we call a preconditioner  {\it feasible} if and only if it is strictly greater than the linear operator under preconditioning, i.e., $\mathbb{M}_{k,\kappa} > T_k$ in Lemma \ref{lemma:preconditioned_DCA:stepsize} or $\mathbb{M} > T$ in Proposition \ref{prop:nopre}.  It is rooted in the requirement for the convergence of the preconditioned DCA with $M_0>0$ as in \eqref{eq:pdca:framework}.  
We would like to emphasize that any finite {\it feasible} preconditioned iterations can still be seen as a whole preconditioned iteration with a modified {\it feasible} preconditioner \cite{DS} (Proposition 3).  
We will design some feasible Jacobi and Richardson preconditioners involving the graph Laplacian $L$ and the normalized graph Laplacian $L_s$. Compared to the symmetric Gauss-Seidel preconditioners, they are more convenient for parallel computing.
Let us begin with the following proposition.

\begin{proposition}\label{prop:diagonal:pre}
	For the graph Laplacian operator $L=D-W$ or the normalized graph Laplacian operator $L_s:=I-D^{-\frac{1}{2}}WD^{-\frac{1}{2}}$,  the following properties are true:
	\begin{equation}
	\begin{cases}
	2D \geq L=D-W, \\
	2I \geq L_s=I-D^{-\frac{1}{2}}WD^{-\frac{1}{2}}.
	\end{cases}
	\end{equation}
\end{proposition}
\begin{proof}
	Since the diagonal matrix $D > 0$, we have
	\[
	2I \geq I-D^{-\frac{1}{2}}WD^{-\frac{1}{2}} \Leftrightarrow 2D^{\frac{1}{2}} ID^{\frac{1}{2}}  \geq  D^{\frac{1}{2}} (I-D^{-\frac{1}{2}}WD^{-\frac{1}{2}}) D^{\frac{1}{2}} \Leftrightarrow 2D \geq L=D-W.
	\]
	Let's focus on proving the case $2D \geq D-W$. It can be readily checked that
	\[
	2D \geq D-W \Leftrightarrow D +W \geq 0.
	\]
	Actually, the positive semidefiniteness of $D+W$ can be deduced from the fact that $D+W$ is diagonal dominant by the definition of the graph Laplacian $L$. 
\end{proof}

For the damped Jacobi or Richardson preconditioner, let us take the operator $T$ in \eqref{eq:dca:linear:system:ori} for example. Suppose the discretization of $T$ is $D - E - E^{*}$ where $D$ is the diagonal part,  $-E$ represents the strictly lower triangular part, and $-E^*$ is the strictly upper triangular part. With a positive definite diagonal matrix $G=\text{diag}(g_{ii})$, the damped Jacobi or Richardson preconditioned iteration can be formulated as follows
\begin{equation}\label{eq:general:diag:pre}
u^{t+1}=u^t + G^{-1}(b-(D-E-E^*)u^t).
\end{equation}
One can write \eqref{eq:general:diag:pre} element-wisely as 
\begin{equation}\label{eq:dampedJacobi:richardson}
u^{t+1}_i = u^{t}_i + \frac{1}{g_{ii}}\left( b_i-\sum_{ j}a_{ij}u_j^{t}\right).
\end{equation}
While $G=D$, \eqref{eq:dampedJacobi:richardson} becomes the classical Jacobi preconditioned iteration and it can also be rewritten as $ u^{t+1}_i =( b_i-\sum_{ j\neq i}a_{ij}u_j^{t})/a_{ii}$. However, the diagonal preconditioner $M=D$ usually did not satisfy the \emph{feasible} condition $D > T$. We thus turn to the damped or \emph{perturbed} Jacobi preconditioned iteration. 

Actually, we can choose the following diagonal preconditioner for \eqref{eq:dca:linear:system:ori}
\begin{equation}\label{eq:pre:unnorm:jacob}
\begin{cases}
G_{Jd} = 2\text{diag}(T)=2(\epsilon D+\eta \Lambda +cI), \quad &\text{damped Jacobi preconditioner}, \\
G_{Jp} = 2 {\epsilon} D+\eta \Lambda +cI + \delta_0 I, \quad &\text{perturbed Jacobi preconditioner},
\end{cases}
\end{equation}
where $\delta_0>0$ is a tiny constant.  Both preconditioners in \eqref{eq:pre:unnorm:jacob} are {\it feasible}, since according to Proposition \ref{prop:diagonal:pre},  it can be checked that 
\[
G_{Jd}-T = 2(\epsilon D+\eta \Lambda +cI) - (\epsilon L+\eta \Lambda +cI)=\epsilon(2D-L) +\eta \Lambda + cI \geq cI >0, 
\]
and
\[
G_{Jp}-T = (2\epsilon D+\eta \Lambda +cI+\delta_0I) - (\epsilon L+\eta \Lambda +cI)=\epsilon(2D-L) +\delta_0 I \geq \delta_0I >0.
\]

We will also consider the following generalized Richardson preconditioner
\begin{equation}\label{eq:pre:unnorm:richard}
G_{R} = \epsilon \lambda_{\max} I + \eta \Lambda + cI + \delta_0I.
\end{equation}
Here, $\delta_0>0$ is as before. The maximum eigenvalue of $L$, i.e., $\lambda_{\max}$, can be estimated by the power method. We further have
\[
G_{R} -T = \epsilon\lambda_{\max} I + \eta \Lambda + cI + \delta_0I -(\epsilon L+\eta \Lambda +cI) = \epsilon (\lambda_{\max} I-L) + \delta_0I \geq \delta_0I >0,
\]
with $(\lambda_{\max} I-L) \geq 0$. Henceforth, we call $G_{R}$ or any other Richardson preconditioner as {\it generalized} Richardson preconditioner since they are not of the form $\alpha I$ with scalar $\alpha>0$ as the usual Richardson preconditioner due to $\Lambda$.

Similarly, for $T = \epsilon L_s + \eta \Lambda +cI $ with  the normalized graph Laplacian $L_s$, we prefer the following \emph{feasible} preconditioners
\begin{equation}\label{eq:pre:norm:jacob:richard}
\begin{cases}
G_{Jd,s} = 2\text{diag}(T)=2(\epsilon I+\eta \Lambda +cI), \quad &\text{damped Jacobi preconditioner}, \\
G_{Jp,s} = 2\epsilon I+\eta \Lambda +cI + \delta_{0} I, \quad &\text{perturbed Jacobi preconditioner}, \\
G_{R,s} = \epsilon \lambda_{s,\max} I + \eta \Lambda + cI + \delta_{0}I, \quad &\text{generalized Richardson preconditioner}.
\end{cases}
\end{equation}
Here,  $\lambda_{s,\max}$ is the maximum eigenvalue of $L_s$.

For the linear operator $T_k$ in \eqref{eq:pdca:framework:stepsize} with the unnormalized graph Laplacian $L$, we choose the following preconditioner
\begin{equation}\label{eq:pre:unnorm:jacob:richard:k}
\begin{cases}
G_{Jd,k} = 2(I + k(\epsilon D+\eta \Lambda +cI)), \quad &\text{damped Jacobi preconditioner}, \\
G_{Jp,k} = I+ k(2\epsilon D +\eta \Lambda +cI) + \delta_0 I, \quad &\text{perturbed Jacobi preconditioner}, \\
G_{R,k} = I +k( \epsilon \lambda_{\max} I + \eta \Lambda + cI) + \delta_0I, \quad &\text{generalized Richardson  preconditioner}.
\end{cases}
\end{equation}
For $T_k$ in \eqref{eq:pdca:framework:stepsize} with the normalized Laplacian $L_s$, we choose 
\begin{equation}\label{eq:pre:norm:jacob:richard:k}
\begin{cases}
G_{Jd,ks} = 2(1+ k(\epsilon +c)I+k\eta \Lambda), \quad &\text{damped Jacobi preconditioner}, \\
G_{Jp,ks} = (2k\epsilon  +1) I+ k(\eta \Lambda +cI) + \delta_{0} I, \quad &\text{perturbed Jacobi preconditioner}, \\
G_{R,ks} = I +k( \epsilon \lambda_{s,\max}  I + \eta \Lambda + cI) + \delta_{0} I, \quad &\text{generalized Richardson preconditioner}.
\end{cases}
\end{equation}

\section{Global Convergence and Kurdyka-\L ojasiewicz (KL) Property}\label{sec:global:conv}
The preconditioners above bring out flexibility.
In this section, we will prove the global convergence of the proposed preconditioned algorithms with any finite feasible preconditioned iterations. By \cite{ABS,DS}, we mainly need to verify the KL property of the original functional $F(u)$.

Now, let us denote  the labeling variables ${u} $, ${y}$  as vectors  $u, y\in \mathbb{R}^{MN}$. 
We will need to use some necessary tools from the standard convex and variational analysis \cite{Roc1}.
The graph of an extended real-valued function $f :\mathbb{R}^n \to \mathbb{R} \cup \{+\infty\}$ is defined by 
\[
\gph f := \{(x,s) \in \mathbb{R}^n\times \mathbb{R} : s = f(x)\}. 
\]
Let $h: \mathbb{R}^n \rightarrow \mathbb{R}\cup \{+\infty\}$ be a proper lower semicontinuous function. Denote $\dom h: = \{ x \in\mathbb{R}^n: \ h(x) < +\infty  \}$. 
For each $x\in\dom f$, 
the limiting-subdifferential of $h$ at $x\in \mathbb{R}^n$, written $\partial f$ are well-defined as in \cite{Roc1} (see Page 301, Definition 8.3(a)). Additionally, if $h$ is continuously differentiable,
then the subdifferential reduces to the usual gradient $\nabla h$, which is our case.
For the global and local convergence analysis, we also need the following  Kurdyka-\L ojasiewicz (KL) property and KL exponent. While the KL properties can help obtain the global convergence of iterative sequences,  the KL exponent can help provide a local convergence rate. 
\begin{definition}[KL property, KL function and KL exponent \cite{ABS}] \label{def:KL}
	A proper closed function $h$ is said to satisfy the KL property at $\bar x \in \dom \partial h$ if there exists $\nu \in (0,+\infty]$, a neighborhood $\mathcal{O}$ of $\bar x$, and a continuous concave function $\psi: [0,\nu) \rightarrow [0,+\infty)$ with $\psi(0)=0$ such that:
	\begin{itemize}
		\item [{\emph{(i)}}] $\psi$ is continuous differentiable on $(0,\nu)$ with $\psi'>0$ over $(0,\nu)$;
		\item [{\emph{(ii)}}] for any $x \in \mathcal{O}$ with $h(\bar x) < h(x) <h(\bar x) + \nu$, one has
		\begin{equation}\label{eq:kl:def}
		\psi'(h(x)-h(\bar x)) \dist(0,\partial h(x)) \geq 1. 
		\end{equation}
	\end{itemize}
	A proper closed function $h$ satisfying the KL property at all points in $\dom \partial h$ is called a KL function. Furthermore, for a proper closed function $h$ satisfying the KL property at $\bar x \in \dom \partial h$,  if $\psi$ in \eqref{eq:kl:def} can be chosen as $\psi(s) = cs^{1-\theta}$ for some $\theta \in [0,1)$ and $c>0$, i.e., there exist  $\bar c, \epsilon >0$ such that
	\begin{equation}\label{eq:KL:exponent:theta:exam}
	\dist(0,\partial h(x)) \geq \bar c (h(x)-h(\bar x))^{\theta}
	\end{equation}
	whenever $\|x -\bar x\| \leq \epsilon$ and $h(\bar x) < h(x) <h(\bar x) + \nu$, then we say that $h$ has the KL property at $\bar x$ with exponent $\theta$.
	If $h$ has the KL property with exponent $\theta$ at any $ \bar x \in \dom \partial h$, we call $h$ a KL function with exponent $\theta$. 
\end{definition}

For KL functions, the following semialgebraic functions can help the KL analysis.
\begin{definition}[Semialgebraic set and Semialgebraic function \cite{ABS}]\label{definition:semialgebraic}
	A subset $S$ of $\mathbb{R}^n$ is called a real semialgebraic set if there exists a finite number of real polynomial functions $P_{i,j}, \ Q_{i,j}: \mathbb{R}^n \rightarrow \mathbb{R}$, such that 
	\[
	S = \bigcup_{i=1}^p \bigcap_{j=1}^q\{ x \in \mathbb{R}^n: P_{i,j} =0, \ Q_{i,j} >0\}.
	\]
	A function $ f:\mathbb{R}^n\to\mathbb{R}\cup\{+\infty\}$ is semialgebraic if its graph is a semialgebraic set of $\mathbb{R}^{n+1}$. 
\end{definition}
A very useful conclusion is that a semialgebraic function has the KL property and is a KL function \cite{ABS} (Lemma 2.3). 
We also need the level boundedness. A function $F: \mathbb{R}^n \rightarrow [-\infty, +\infty]$ is level-bounded (see Definition 1.8 \cite{Roc1}) if lev$_{\leq \alpha}F: =\{ u: F(u) \leq \alpha\}$ is bounded (or possibly empty). The level boundedness is usually introduced for discussions of the existence of the minimizer. It is weaker than coerciveness. For example, the indicator function $I_{[0,1]}(x)$ in $\mathbb{R}$ is level bound but is not coercive.

We now turn to the energy decreasing of the proposed preconditioned DCA \eqref{eq:pdca:framework} before the discussion of its global convergence. Inspired by \cite{DS} (the proof of Theorem 2) and \cite{LPH} (Theorem 3), we have the following theorem.
\begin{theorem}\label{thm:decay:energy}
	For the preconditioned DCA \eqref{eq:pdca:framework} with convex $P_1$ and $P_2$, we have 
	the following  energy  decay: $F(u^{t+1}) \leq F(u^t)$. Furthermore, while $u^{t+1} \neq u^t$, we have $F(u^{t+1}) < F(u^t)$.
\end{theorem}
\begin{proof}
	For $P_1(u)$ and $P_2(u)$ in \eqref{eq:pdca:framework} with positive definite $M_0$, for $\xi^t \in \partial P_2(u)|_{u=u^t}$, by the definition of $u^{t+1}$ in \eqref{eq:pdca:framework}, we have
	\begin{align}
	&P_1(u^{t+1}) -\langle \xi^t, u^{t+1}\rangle + \frac{1}{2}\|u^{t+1}-u^t\|_{M_0}^2 \leq P_1(u^t)  -\langle \xi^t, u^{t}\rangle \notag \\
	&\Leftrightarrow P_1(u^{t+1}) + \frac{1}{2}\|u^{t+1}-u^t\|_{M_0}^2 \leq P_1(u^t) + \langle \xi^t, u^{t+1}-u^t\rangle. \label{eq:energy:decay1}
	\end{align}
	By the convexity of $P_2$, we have 
	\begin{equation}\label{eq:sub:energy}
	\langle \xi^t, u^{t+1}-u^t\rangle \leq P_2(u^{t+1})-P_2(u^t).
	\end{equation}
	Substituting \eqref{eq:sub:energy} into \eqref{eq:energy:decay1}, we have
	\begin{align*}
	&P_1(u^{t+1}) -P_2(u^{t+1}) +\frac{1}{2}\|u^{t+1}-u^t\|_{M_0}^2 \leq  P_1(u^t) -P_2(u^t) \\
	&\Leftrightarrow F(u^{t+1}) + \frac{1}{2}\|u^{t+1}-u^t\|_{M_0}^2 \leq F(u^t).
	\end{align*}
	We then get this theorem with positive definiteness of $M_0$.
\end{proof}


We will also consider the convergence rate of the sequence $\{u^t\}$ if the KL exponent of $F(u)$ is known. This kind of convergence rate analysis is standard; see \cite{ABS} for a more comprehensive analysis. We give the following theorem \cite{ABS}   and the proof is omitted here.
\begin{theorem}[convergence rate]\label{thm:rate}
	Assume that $u^t$ converges to $\bar{u}$ and $F(u)$ has the KL property at $\bar{u}$ with $\psi(s)=cs^{1-\theta},\theta\in\left[0,1\right)$, $c>0$. Then the following estimations hold:
	\begin{enumerate}
		\item If $\theta=0$ then the sequence $u^t$ converges in a finite number of steps,
		\item If $\theta\in\left(0,\frac12\right]$ then there exist $c>0$ and $\tau\in\left[0,1\right)$, such that $\|u^t-\bar{u}\|_2\le c\tau^t$,
		\item If $\theta\in\left(\frac12,1\right]$ then there exists $c>0$, such that $\|u^t-\bar{u}\|_2\le c t^{-\frac{1-\theta}{2\theta-1}}$.
	\end{enumerate}
\end{theorem}
With these preparations, we finally arrive at the following theorem for the global convergence and local convergence rate of the preconditioned DCA \eqref{eq:pdca:framework}.
\begin{theorem}\label{thm:energy:decrease}
	The energy functional $F(u)$ in \eqref{eq:discrete:fu} is level bounded and is a KL function with KL exponent $1-R(n,d)$ with $R(n,d) =  1/(3^{3(MN-1)}4)$. If $\bar{u}$ is the isolated zero of $F(u)$, its KL exponent is $1-\chi(n,d)$ around $\bar{u}$ with $\chi(n,d) = 1/(3^{MN}+1)$. The iteration sequence $\{u^t\}$ is bounded and globally convergent. The local convergence rate of preconditioned DCA \eqref{eq:pdca:framework}
	is sublinear. Here we assume $u \in \mathbb{R}^{MN}$ with $n=MN$ and $d=4$ being the polynomial order of $W(u)$.
\end{theorem}
\begin{proof}
	The level boundedness of $F(u)$ follows from the coercivity of $W(u)$, the positive semidefiniteness of the graph Laplacian $L$, and the weighted matrix $\Lambda$. It can be deduced by the contradiction argument. Since for any positive constant $\alpha>0$, if there exists a sequence $\{u^{t_i}\}$ such that $F(u^{t_i})\leq \alpha$ and $\|u^{t_i}\|_{2} \rightarrow +\infty$, by the coercivity of $W(u)$, we see that $F(u)\geq \mathbb{W}(u) \rightarrow +\infty$, which is contradicted $F(u)\leq \alpha$.
	Together with the energy decreasing of $F(u^t)$ in Theorem \ref{thm:energy:decrease} and the level boundedness of $F$, we conclude that the sequence $\{ u^t\}$ is bounded.
	
	Furthermore,  the graph of $F$, i.e., $\gph F := (u,z)$ satisfies
	\begin{align}
	&z = \sum_{ij}\frac{\epsilon}{2}w_{ij}(u(i)-u(j))^2 + \frac{1}{4\epsilon}\sum_i (u(i)^2-1)^2 + \frac{\eta}{2}\sum_i\Lambda(i)(u(i)-y(i))^2 \\
	&\Longleftrightarrow z - \sum_{ij}\frac{\epsilon}{2}w_{ij}(u(i)-u(j))^2 + \frac{1}{4\epsilon}\sum_i (u(i)^2-1)^2 + \frac{\eta}{2}\sum_i\Lambda(i)(u(i)-y(i))^2=0.
	\end{align}
	The graph $\gph F$ is a semialgebraic set since
	\[
	\left\{(u,z) \ | \ z - \sum_{ij}\frac{\epsilon}{2}w_{ij}(u(i)-u(j))^2 + \frac{1}{4\epsilon}\sum_i (u(i)^2-1)^2 + \frac{\eta}{2}\sum_i\Lambda(i)(u(i)-y(i))^2=0\right\}.
	\]
	We thus conclude that $F(u)$ is a semialgebraic function, which is a KL function.  Together with boundedness of $\{ u^t\}$,  we see $\{ u^t\}$ is globally convergent \cite{DS} (Theorem 2 and Theorem 3).
	
	Moreover, since $F(u)$ is a polynomial with order $d=4$ and $u \in \mathbb{R}^{MN}$ with $MN$ variables, we conclude that the KL exponent of $F$ is $1-R(n,d)$ by \cite{LMP} (Lemma 2). While $\bar u$ is the isolated zeros of $F$,  the KL exponent of $F$ in the neighborhood of $\bar{u}$ is $1-\chi(n,d)$ still by \cite{LMP} (Lemma 2).  
	Since $1-R(n,d) \in (1/2,1]$ and $1-\chi(n,d) \in (1/2,1]$, according to Theorem \ref{thm:rate}, the local convergence rate is sublinear. 
\end{proof}

Now, let us end this section with the following remark before showing the efficiency of \eqref{eq:pdca:framework} with numerical tests.
\begin{remark}
	Let's give a comparison to some existing convergence analyses of convex splitting methods \cite{Ey, LB} for the convergence of DCA with arbitrary step size $k>0$ and no preconditioning \eqref{eq:dca:stepsize}. The convergence analysis based on KL properties here only assumes the convexity of $P_2(u)$ and does not need the strictly uniform convexity of $P_2(u)$ or the semi-convexity of $F(u)$ as required in \cite{Ey}. The convergence is also guaranteed with unconditional stability compared to conditional stability as in \cite{LB}, which is based on the semi-implicit splitting method.   
\end{remark}
\section{Numerical results}\label{sec:num}

This section is divided into two parts.  The first part focuses on image segmentation with GPU including discussions on the search window. It can help construct different graph Laplacian matrices. The second part focuses on data clustering with \eqref{eq:pdca:framework}.

\subsection{Image segmentation by GPU}





\subsubsection{Comparing with SVD}
\begin{figure*}[!htbp]
	\centering
	\subfloat[two cow: 320$\times$312 \label{fig:subfig:Originalcow}]{
		\includegraphics[width=1.8in, origin=br]{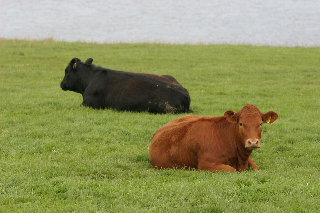}}
	\subfloat[DICE=0.9 (svd1)\label{fig:subfig:Cow__SVD_90per_295}]{
		\includegraphics[width=1.8in, origin=br]{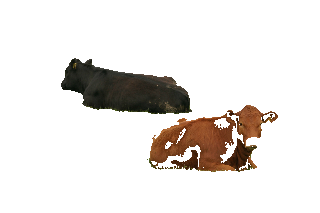}}
	\subfloat[svd1 final result \label{fig:subfig:Cow__SVD_96per_366}]{
		\includegraphics[width=1.8in, origin=br]{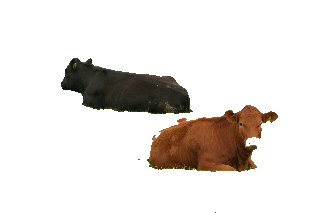}}
	\\
	\subfloat[hand-labeled image\label{fig:subfig:Cowprior}]{
		\includegraphics[width=1.8in, origin=br]{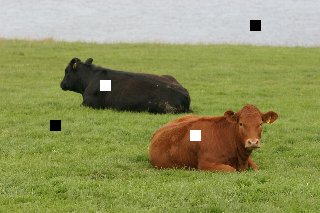}}
	\subfloat[DICE=0.9 (pDCA)\label{fig:subfig:Cow_pDCA_90per_79} ]{
		\includegraphics[width=1.8in, origin=br]{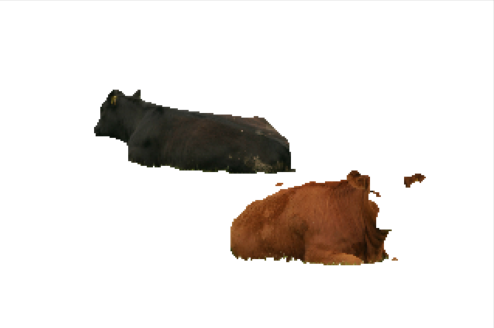}}
	\subfloat[pDCA final result\label{fig:subfig:Cow_pDCA_97per_404}]{
		\includegraphics[width=1.8in, origin=br]{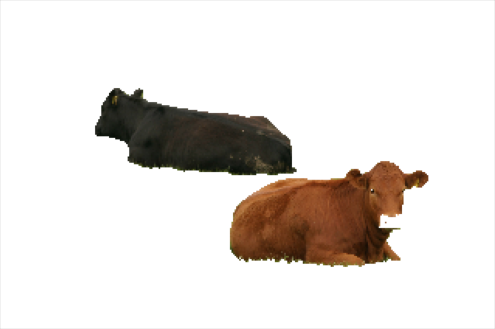}}
	\caption{Image {\bf{a}} is the original two cows image. Image {\bf{d}} contains the prior information.  Images {\bf{b}} and  {\bf{c}} are the segmented images obtained by SVD-Nystr\"{o}m based method (svd1) with DICE=0.9 and stable DICE correspondingly.  Images {\bf{e}} and  {\bf{f}} are the segmented images obtained by the proposed pDCA with DICE=0.9 and stable DICE correspondingly.  }
	\label{fig:Segmentationbygeneralwindow_cow} 
\end{figure*}

\begin{figure*}[!htbp]
	\centering
	\subfloat[stone: 640$\times480$ \label{fig:subfig:Originalstone}]{
		\includegraphics[width=1.8in, origin=br]{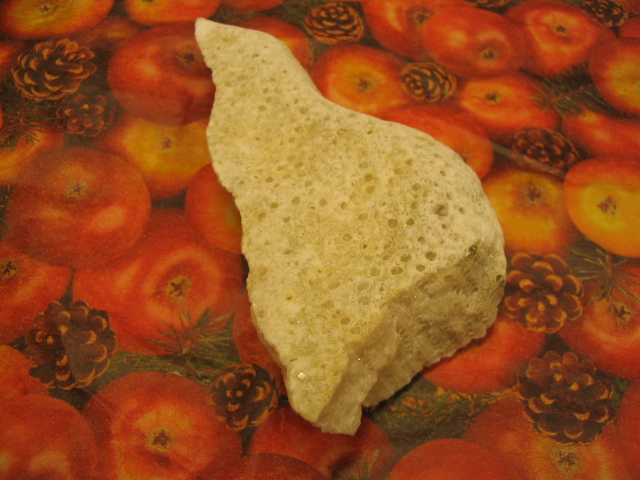}}
	\subfloat[DICE=0.9 (svd1)\label{fig:subfig:Stone_SVD_90per_89}]{
		\includegraphics[width=1.8in, origin=br]{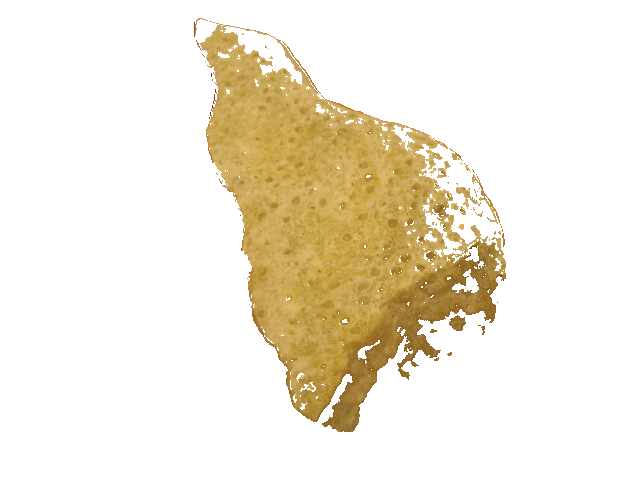}}
	\subfloat[svd1 final result \label{fig:subfig:Stone_SVD_96per_112}]{
		\includegraphics[width=1.8in, origin=br]{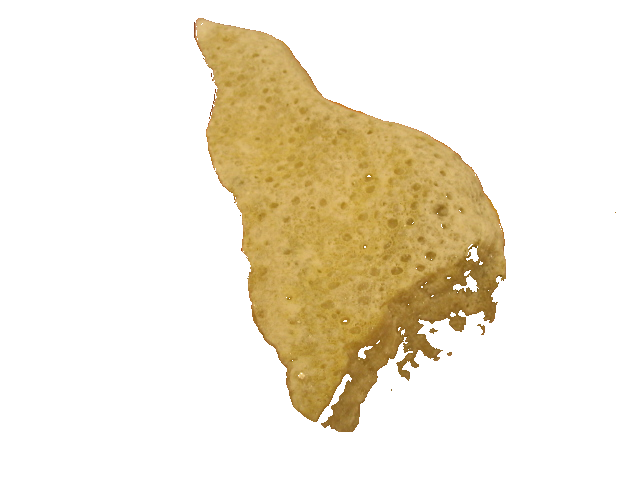}}
	\\
	\subfloat[hand-labeled image\label{fig:subfig:Stoneprior}]{
		\includegraphics[width=1.8in, origin=br]{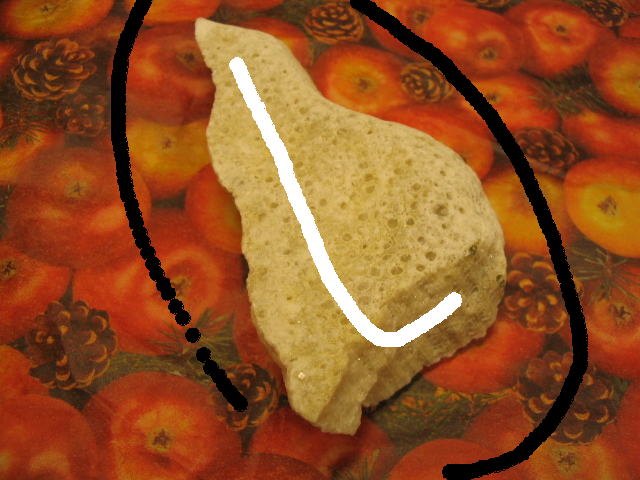}}
	\subfloat[DICE=0.9 (pDCA)\label{fig:subfig:Stone_pDCA_90per_126} ]{
		\includegraphics[width=1.8in, origin=br]{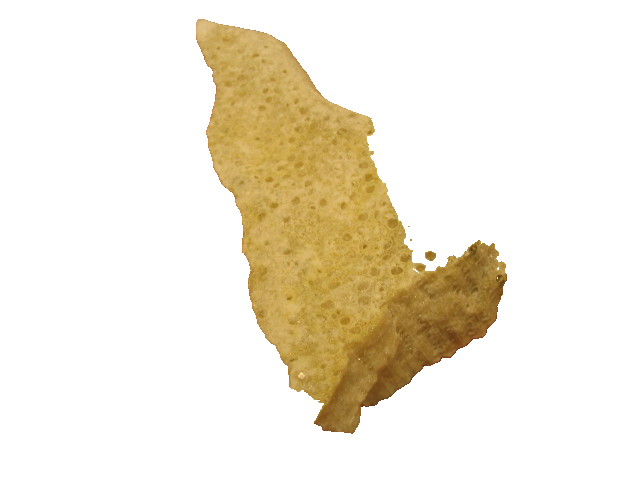}}
	\subfloat[pDCA final result\label{fig:subfig:Stone_pDCA_99per_292}]{
		\includegraphics[width=1.8in, origin=br]{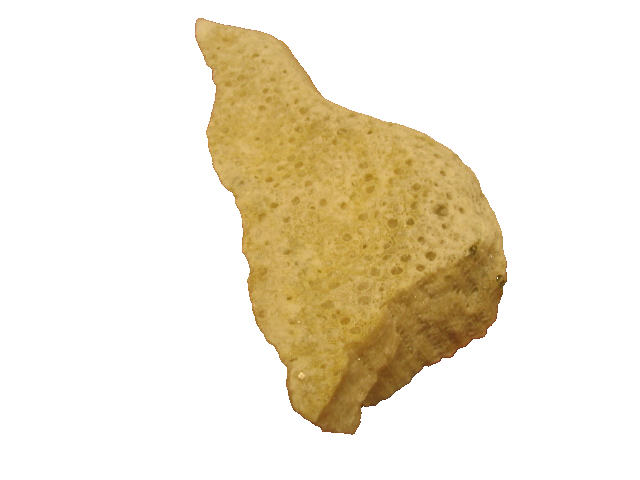}}
	\caption{Image {\bf{a}} is the original stone image. Image {\bf{d}} contains the prior information.  Images {\bf{b}} and  {\bf{c}} are the segmented images obtained by SVD-Nystr\"{o}m based method (svd1) with DICE=0.9 and stable DICE correspondingly.  Images {\bf{e}} and  {\bf{f}} are the segmented images obtained by the proposed pDCA with DICE=0.9 and stable DICE correspondingly. }
	\label{fig:Segmentationbygeneralwindow_stone} 
\end{figure*}

\begin{figure*}[!htbp]
	\centering
	\subfloat[clover: 1000$\times$800 \label{fig:subfig:Originalclover}]{
		\includegraphics[width = 1.3in, origin=br]{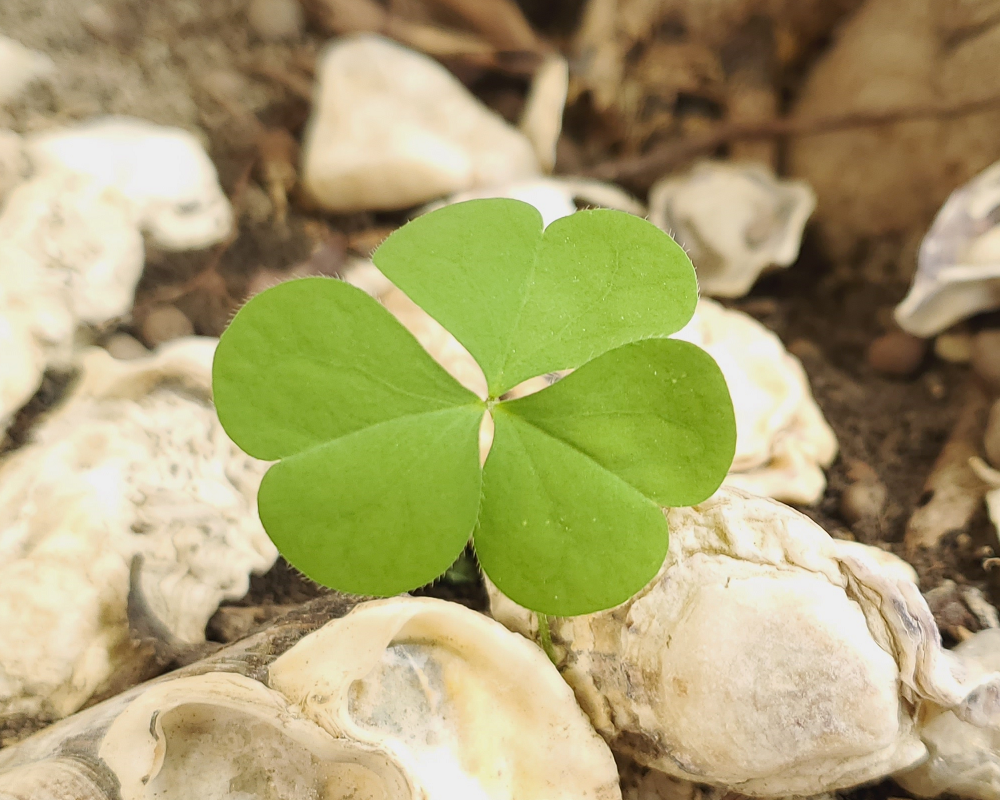}}
	\subfloat[hand-labeled image\label{fig:subfig:Cloverprior}]{
		\includegraphics[width=1.3in, origin=br]{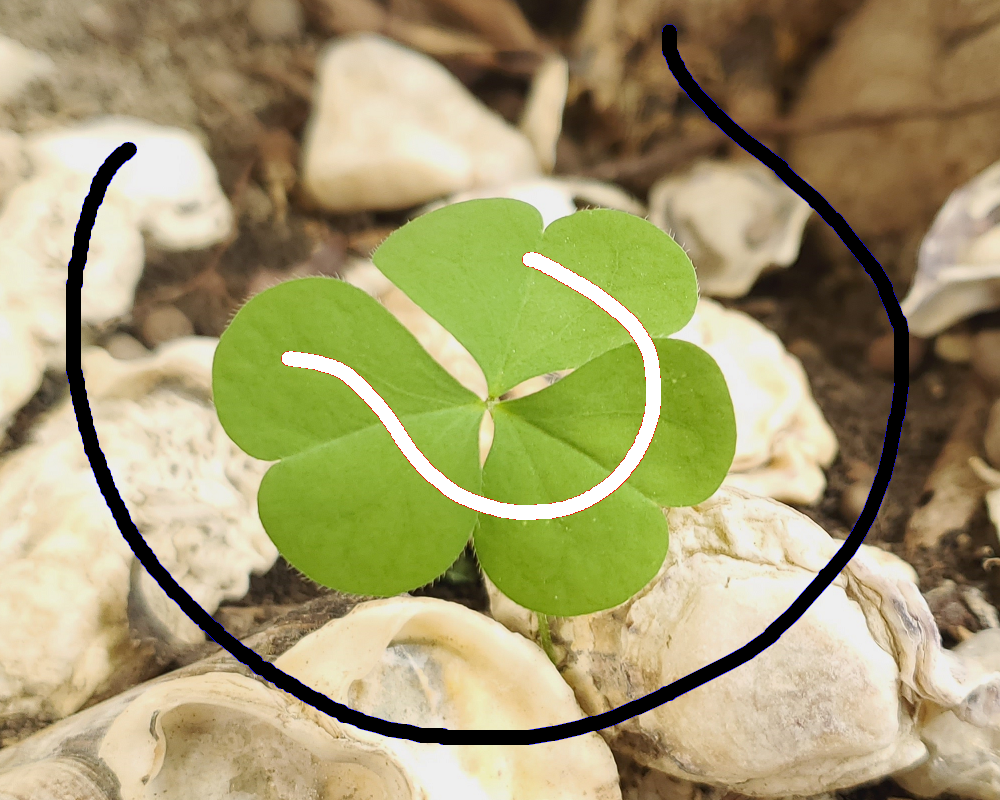}}
	\subfloat[DICE=0.9\label{fig:subfig:Clover_pDCA_90per_333}]{
		\includegraphics[width=1.3in, origin=br]{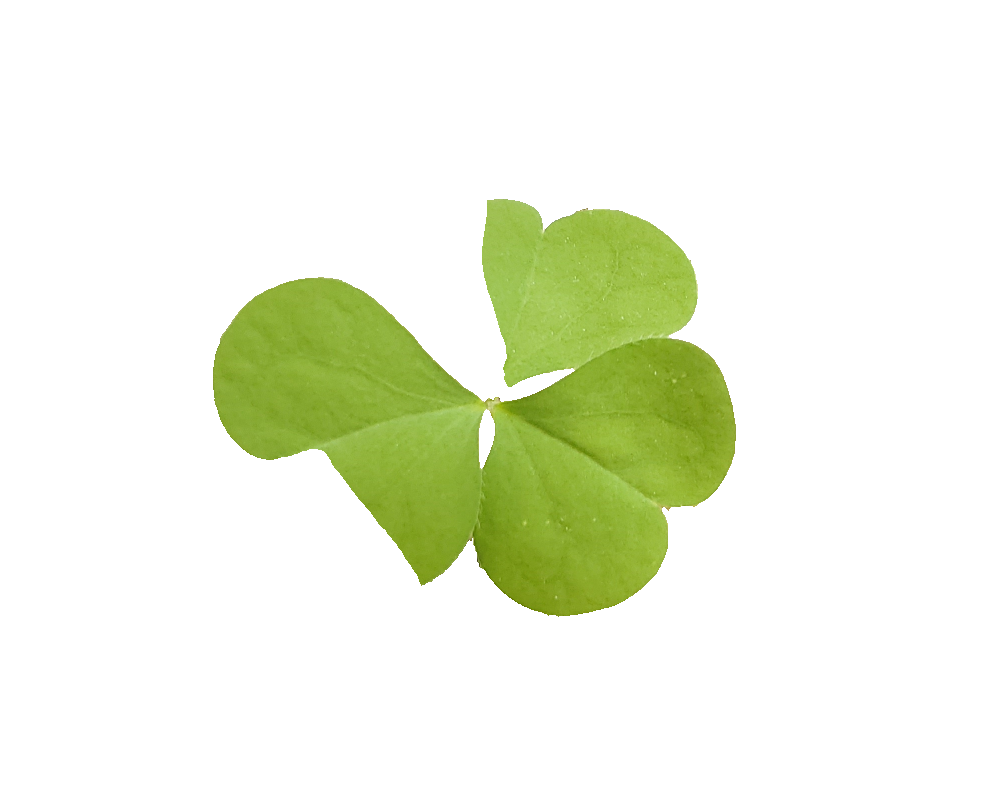}}
	\subfloat[pDCA final result ]{
		\includegraphics[width=1.3in, origin=br]{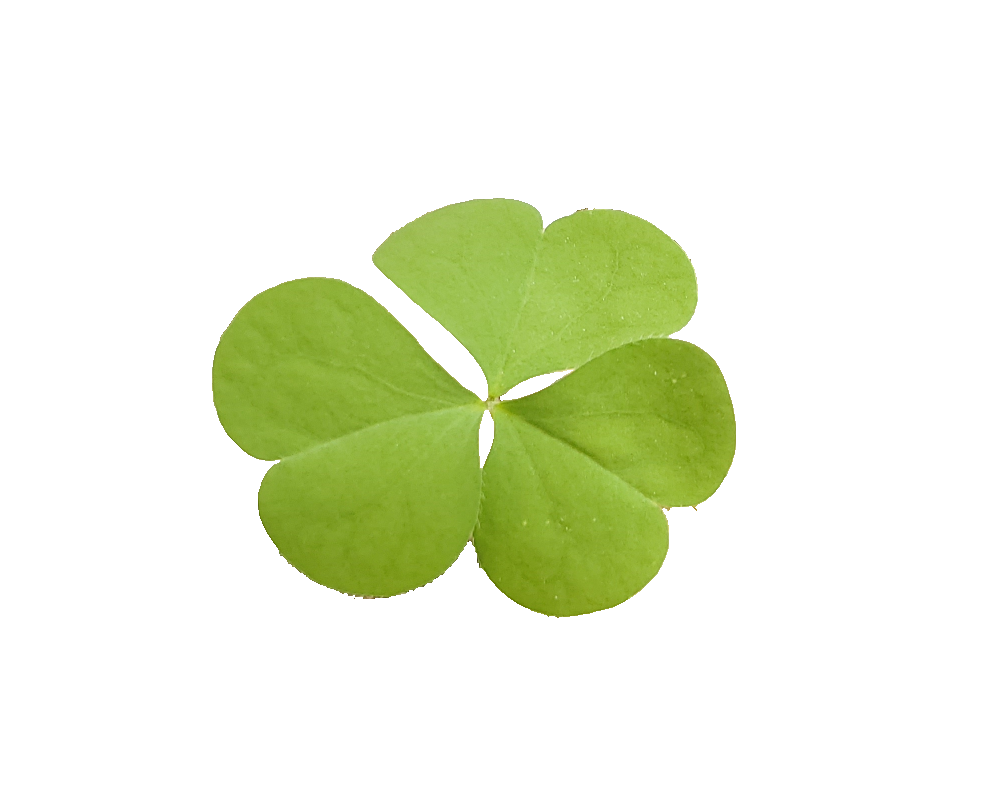}}
	
	\caption{Image {\bf{a}} is the original clover image. Image {\bf{b}} contains the prior information.  Images {\bf{c}} and  {\bf{d}} are the segmented images obtained by the proposed pDCA with DICE=0.9 and stable DICE correspondingly.}
	\label{fig:Segmentationbygeneralwindow_clover} 
\end{figure*}

In comparison with SVD-Nystr\"{o}m method \cite{BF,GMBFP,LB},  we provide several examples to show our advantages in segmentation. In the experiments, the SVD-Nystr\"{o}m algorithm and the proposed \eqref{eq:pdca:framework} are executed on  the workstation 
(CPU: Intel(R) Xeon(R) CPU E5-2650 v4@2.20GHz; GPU: NVIDIA Corporation GM204GL [Quadro M4000], 8GB). There are different sizes of images including the smallest image of two cows with size 320$\times$312 in Figure \ref{fig:subfig:Originalcow},  the largest clover image with size 1000$\times$800 in Figure \ref{fig:subfig:Originalclover}, and the stone image with size 640$\times$480 in Figure \ref{fig:subfig:Originalstone}. Regarding the accuracy of segmentation, we take the following DICE and  Jaccard coefficients to show the performance between the two algorithms
\begin{equation*}
\text{DICE} = \frac{2\text{TP}}{\text{FP}+2\text{TP}+\text{FN}},\qquad
\text{Jaccard} = \frac{\text{TP}}{\text{FP}+\text{TP}+\text{FN}}.
\end{equation*}
Here, `TP', `FP', and `FN' respectively present the pixels in both ground truth and segmentation, the pixels in segmentation but not in the ground truth, and the pixels in ground truth but not in segmentation. 
The algorithms along with the corresponding parameters are as follows.
\begin{enumerate}
	\item svd1: the SVD-Nystr\"{o}m based convex splitting method \cite{GMBFP} with $\epsilon = 100$, $c = 11 $, $\eta = 100$, $dt = 0.001$. 
	\item svd2: the SVD-Nystr\"{o}m based difference of convex functions algorithm as in Remark  \ref{rem:dca:constantcoefficient} with $\epsilon = 100$, $c = 120 $, $\eta = 100$. svd2 did not have a step size or can be seen as having an infinite step size.
	\item pDCA: $\epsilon = 100$, $c = 11 $, $\eta = 100$, 4 preconditioned Richardson iterations are chosen,  $150$ iterations of power method for the largest eigenvalue of normalized Laplacian.
\end{enumerate}

\begin{remark}\label{rem:code:souce}
	For the parameter settings in image segmentation experiments, we mainly refer to the publicly available  SVD-Nystr\"{o}m based  Matlab code \cite{BF, GMBFP} \footnote{\url{https://users.math.msu.edu/users/merkurje/GL_segmentation.zip}, accessed on
		December 2021.}. We choose the same parameter $\epsilon$ as the code. For other parameters including $\eta$, we choose them according to our numerical experiences together with code of \cite{BF, GMBFP}.
\end{remark}

\begin{figure*}[!htbp]
	\centering
	\subfloat[red flowers: 481$\times$321 \label{fig:subfig:Originalredflower}]{
		\includegraphics[width=1.8in, origin=br]{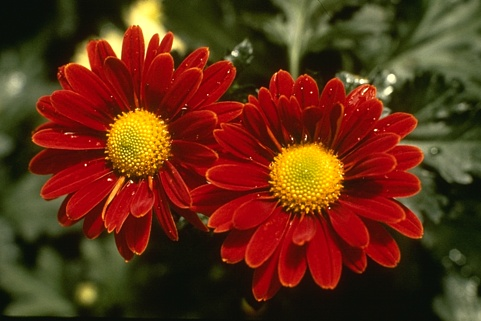}}
	\subfloat[DICE=0.9 (svd1)\label{fig:subfig:Redflower_SVD_90per_494}]{
		\includegraphics[width=1.8in, origin=br]{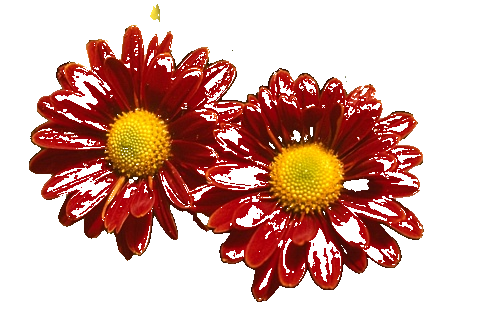}}
	\subfloat[svd1 final result \label{fig:subfig:Redflower_SVD_98per_54}]{
		\includegraphics[width=1.8in, origin=br]{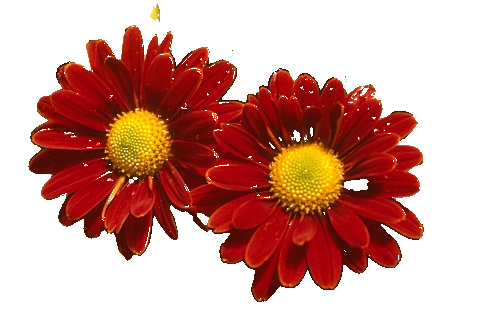}}
	\\
	\subfloat[hand-labeled image\label{fig:subfig:Redflowerprior}]{
		\includegraphics[width=1.8in, origin=br]{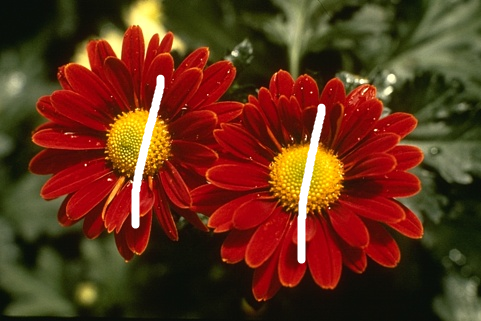}}
	\subfloat[DICE=0.9 (pDCA)\label{fig:subfig:redflower_pDCA_90per_262} ]{
		\includegraphics[width=1.8in, origin=br]{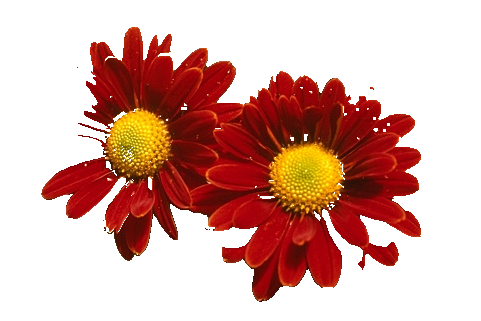}}
	\subfloat[pDCA final result\label{fig:subfig:Redflower_pDCA_98per_788}]{
		\includegraphics[width=1.8in, origin=br]{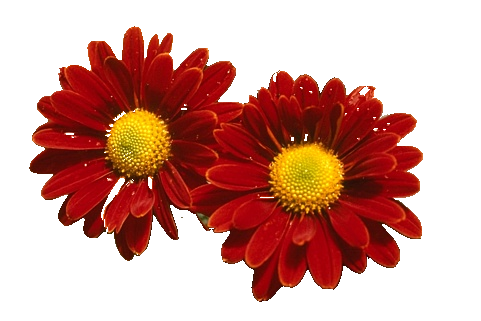}}
	\caption{Image {\bf{a}} is the original red flowers image. Image {\bf{d}} contains the prior information with only labeling on the flower object.  Images {\bf{b}} and  {\bf{c}} are the segmented images obtained by SVD-Nystr\"{o}m based method (svd1) with DICE=0.9 and stable DICE correspondingly.  Images {\bf{e}} and  {\bf{f}} are the segmented images obtained by the proposed pDCA with DICE=0.9 and stable DICE correspondingly.}
	\label{fig:Segmentationbygeneralwindow_redflower} 
\end{figure*}

All compared algorithms including svd1, svd2, and pDCA are semi-supervised. 
As shown in Figures \ref{fig:Segmentationbygeneralwindow_cow}-\ref{fig:Segmentationbygeneralwindow_clover}, there are some priors (or user interactions) in the image, which can ensure the foreground and background of the image in order to get the segmentation. 
Here we use the `square seed' prior for the image `two cows', which is the same as in \cite{LB}. For all other images,  we use `seed' prior such as subfigure (d) of Figures \ref{fig:Segmentationbygeneralwindow_redflower} to \ref{fig:Segmentationbygeneralwindow_clover}    due to their large sizes.

First, we compare the performances of two kinds of algorithms in low sampling or a small searching window in Table \ref{tab:comparison with less data}. Searching windows is the essential factor in constructing the Laplacian matrices for different images with pDCA. In addition, we also include the running time in different stages of the SVD-Nystr\"{o}m based methods, i.e., svd1 and svd2. It can be seen from Table \ref{tab:comparison with less data} that calculating eigenvalues and eigenvectors by SVD with Nystr\"{o}m sampling costs the main computation efforts for svd1 or svd2. Moreover, as shown in Table \ref{tab:comparison with less data}, the proposed svd2 here does not have many advantages over svd1  and we will only compare with the svd1 henceforth.
\begin{table}[h]
	\centering
	\begin{tabular}{|l|ccc|ccc|cc|}\hline
		image & \multicolumn{3}{c|}{two cows}  & \multicolumn{3}{c|}{stone} & \multicolumn{2}{c|}{clover}\\ 
		&svd1      &svd2     &pDCA       &svd1        &svd2        &pDCA       &svd      &pDCA\\\hline
		\multicolumn{9}{|c|}{DICE = 0.90} \\ \hline
		{time}      &74(48)    &52(47)   &\bf{10}         &1193(1137)  &1165(1153)  &\bf{44}         &--            &106\\
		{iteration} &293       &37       &348        &89          &12          &393           &--          &333\\ \hline 
		\multicolumn{9}{|c|}{Stable DICE} \\ \hline
		{time}             &76(48)    &53(47)   &\bf{36}         &1230(1137)  &1168(1153)  &\bf{78}         &--          &228\\
		{iteration}         &313       &39       &1752       &112         &15          &849           &--          &943\\
		{DICE}            &0.9600   &0.9602  &\bf{0.9630}    &0.9600     &0.9611     &\bf{0.9861}    &--           &0.9900\\
		{Jaccard}  &0.9233   &0.9237  &\bf{0.9286}    &0.9230     &0.9232     &\bf{0.9872}    &--          &0.9940\\\hline
	\end{tabular}
	\caption{We compare the running time measured in seconds and the number of iterations between SVD-Nystr\"{o}m based methods (svd1 and svd2) and pDCA. For svd1 and svd2, the sampling parameters are 0.3\% of the size of pictures(`two cows': 200, `stone': 920), and the step sizes are $0.001$ and $\infty$ respectively. Besides, for svd1 and svd2, we also record the time used for SVD (the running time in brackets), which is an important part of the whole program. For pDCA, we choose the smaller searching window which contains $15\times 15$ pixels to construct the graph Laplacian matrix. Since they are too slow, there is no data for either svd1 or svd2 for the large `clover' image. The stable `DICE' is obtained under the principle that the parameter `DICE' does not change or vary in a small range, generally $ 10^{-5}$ within 10 consecutive steps.
	}
	\label{tab:comparison with less data}
\end{table}

Second, when increasing the sampling or magnifying the searching window as in Table \ref{tab:comparison with more data} with figures two cows, stone, and red flowers (size: 481$\times$321), both svd1 and pDCA algorithms will spend more time processing data. However,  compared with SVD-Nystr\"{o}m based method,  pDCA is much faster with the Richardson preconditioner in favor of parallel computing. We can update the label of each pixel on the different threads of the GPU simultaneously. This can help deal with massive computations in a short time without SVD.  

Third, as shown in Figures \ref{fig:Segmentationbygeneralwindow_cow}, \ref{fig:Segmentationbygeneralwindow_redflower}, \ref{fig:Segmentationbygeneralwindow_stone}, and \ref{fig:Segmentationbygeneralwindow_clover}, the proposed pDCA can obtain high-quality segmentation efficiently compared to the current SVD-Nystr\"{o}m based method (svd1) \cite{BF,LB}. 
Also from Tables \ref{tab:comparison with less data} and \ref{tab:comparison with more data}, we see that the advantage of pDCA algorithm is that it remains highly efficient to reach rough segmentations for large images, and also gets a better result than SVD-Nystr\"{o}m based method in the end. In addition, our pDCA algorithm reaches higher accuracy rates in DICE and can get better visual segmentation.

\begin{figure*}[!htbp]
	\centering
	\subfloat[Original image \label{fig:subfig:fullconnectionoriginalimage}]{
		\includegraphics[width=1.3in, origin=br]{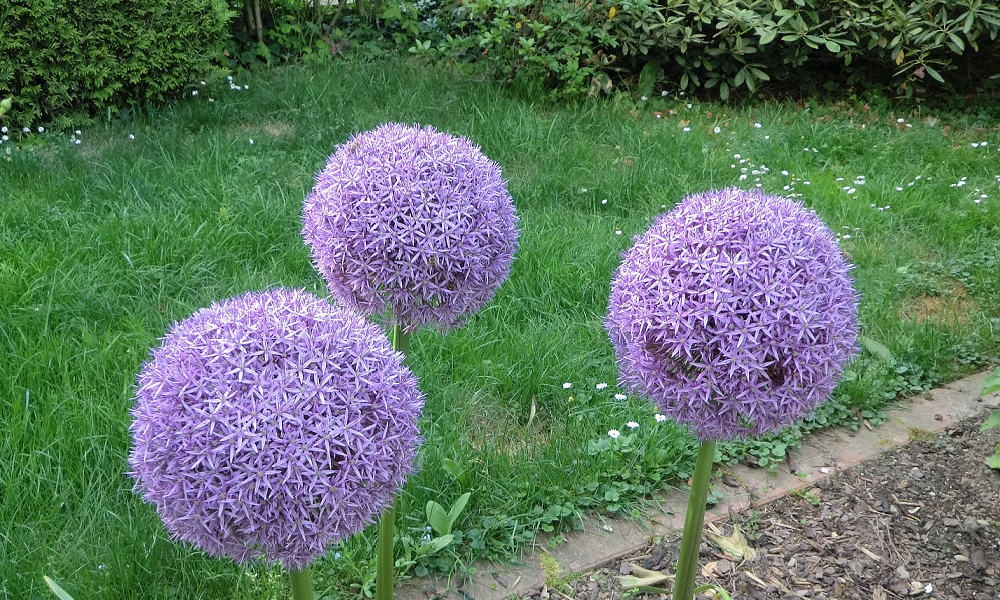}}
	\subfloat[Resize image\label{fig:subfig:resizeimage}]{
		\includegraphics[width=1.3in, origin=br]{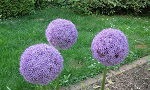}}
	\subfloat[Prior 1\label{fig:subfig:purpleprior1} ]{
		\includegraphics[width=1.3in, origin=br]{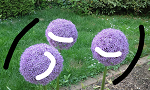}}
	\subfloat[Result 1\label{fig:subfig:fullconnectionruslt1}]{
		\includegraphics[width=1.3in, origin=br]{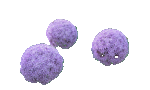}}\\
	\subfloat[Prior 2\label{fig:subfig:purpleprior2} ]{
		\includegraphics[width=1.3in, origin=br]{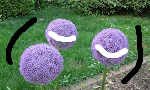}}
	\subfloat[Result 2\label{fig:subfig:fullconnectionruslt2}]{
		\includegraphics[width=1.3in, origin=br]{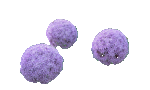}}
	\subfloat[Prior 3\label{fig:subfig:purpleprior3} ]{
		\includegraphics[width=1.3in, origin=br]{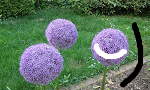}}
	\subfloat[Result 3\label{fig:subfig:fullconnectionruslt3}]{
		\includegraphics[width=1.3in, origin=br]{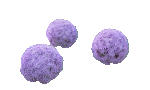}}
	\caption{Subfigure {\bf{a}} is the original 1000$\times$600 purple flower image and we resize it by shrinking it to 15 percent of its original size, i.e., $ 150\times 90$ pixels. Subfigures {\bf{c}}, {\bf{e}}, and {\bf{g}} are different priors of the resized images and subfigures {\bf{d}}, {\bf{f}}, and {\bf{h}} are the corresponding segmented images by pDCA based on different priors. }
	\label{fig:fullconnectionsegmentation} 
\end{figure*}

We now turn to the choices of search windows. A search window that is large enough not only can help obtain more accurate segmentation but also can reduce the dependence on the prior information. However, large search windows also cost massive storage and expensive computations.
In Figure \ref{fig:fullconnectionsegmentation}, we use the entire image as the fully connected search window. It can be seen that the final segmented image does not sensitive to the priors. 
%
\begin{table}[h]
	\centering
	\begin{tabular}{|l|cc|cc|cc|}\hline
		image name & \multicolumn{2}{c|}{two cow} & \multicolumn{2}{c|}{red flowers} & \multicolumn{2}{c|}{stone} \\
		&svd1        &pDCA    &svd1      &pDCA    &svd1     &pDCA\\ \hline
		\multicolumn{7}{|c|}{DICE = 0.90} \\ \hline                 
		time    &255        &\bf{27}      &1621     &\bf{45}      &2988    &\bf{40}\\
		iteration &308        &\bf{79}      &497      &\bf{262}     &88      &126\\ \hline
		\multicolumn{7}{|c|}{Stable DICE} \\ \hline
		time              &265        &\bf{62}      &1727     &\bf{110}     &3425    &\bf{61}\\
		iteration         &367        &403          &590      &903          &192     &292\\
		DICE              &0.9601    &\bf{0.9700} &0.9800  &\bf{0.9830}     &0.9650 &\bf{0.9930}\\
		Jaccard            &0.9234    &\bf{0.9417} &0.9607  &\bf{0.9666} &0.9323 &\bf{0.9861}\\ \hline
	\end{tabular}
	\caption{We compare the running time measured in seconds and the number of iterations between svd1 and pDCA. Different from Table \ref{tab:comparison with less data}, we choose the bigger sampling parameter for sv1 and the larger searching window for pDCA. 
		For svd1, the step size is $0.001$.
		The sampling parameters are 1\% of the size of pictures (`two cows': 680, `red flowers': 1544, `stone': 3072). The size of searching windows depends on the memory of our GPU (`two cows': $55\times 55$, `red flowers': $35\times 35$, `stone': $25\times 25$). The running time for svd1 or pDCA is the running time for the whole program respectively. The stable `DICE' is obtained under the principle that the parameter `DICE' does not change or vary in a small range, generally $ 10^{-5}$ in 10 consecutive steps.}
	\label{tab:comparison with more data}
\end{table}

We give the following remark for implementations to end this subsection. 
\begin{remark}
	We employ the svd1 (or svd2) as in the publicly available implementation of \cite{BF, GMBFP},   which is implemented by Matlab with extremely efficient SVD, while the essential steps including computing weights and preconditioned iterations of pDCA are implemented by CUDA. 
\end{remark}

\subsubsection{Different methods to construct the graph Laplacian matrix for pDCA}
\begin{figure*} [htbp!]
	\centering
	\subfloat[sparse method 1\label{fig:search_window1}]{
		\includegraphics[width=1.5in,  origin=br]{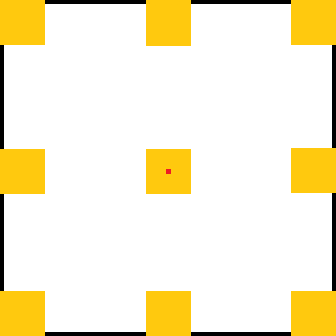}}\ \qquad \qquad
	\subfloat[sparse method 2\label{fig:search_window2}]{
		\includegraphics[width=1.5in,  origin=br]{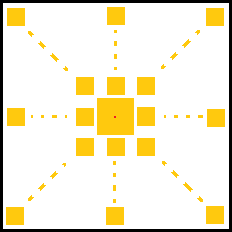}}\ \ 
	\caption{Two different sparse methods have different advantages. Figure \ref{fig:search_window1} uses nine blocks located in nine different directions. Figure \ref{fig:search_window2} is based on figure \ref{fig:search_window1} and reduces the size of blocks except for the central blocks, and adds more same small blocks in the same directions.} 
	\label{fig:search_window} 
\end{figure*}

We now focus on the specially designed sparse search windows for pDCA. Note that svd1 and svd2 only depend on the sampling rate. For nonlocal models on large-size images, it is usually challenging for solving the graph Laplacian matrix $L$. The matrix $L$ is usually too large for storing in computer memory but also costs too much time to compute the action on $u$, i.e., $Lu$. In this section, we propose some specially designed search windows to make $L$ more compact as in Figure \ref{fig:search_window}. These search windows are complementary to the usual box search window as in Figure \ref{patch_and_window}b. They can capture global information for images with disconnected components while keeping sparse structures convenient for storage and computations.

\begin{figure*}[!htb]
	\centering
	\subfloat[Image: 400$\times$300 \label{fig:subfig:Originalyellowflower}]{
		\includegraphics[width=1.8in, origin=br]{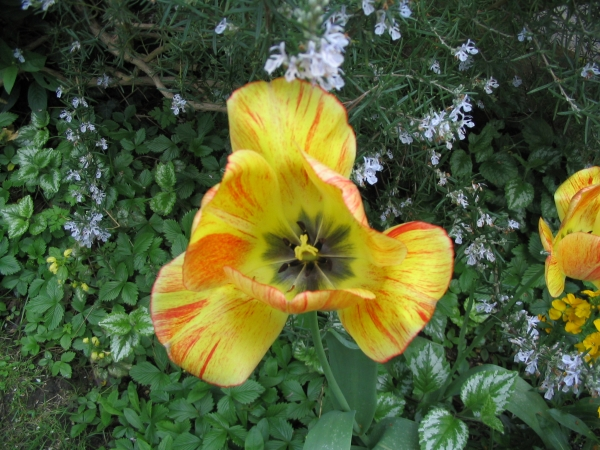}}
	\subfloat[Ground truth\label{fig:subfig:Yellowflowergroundtruth}]{
		\includegraphics[width=1.8in, origin=br]{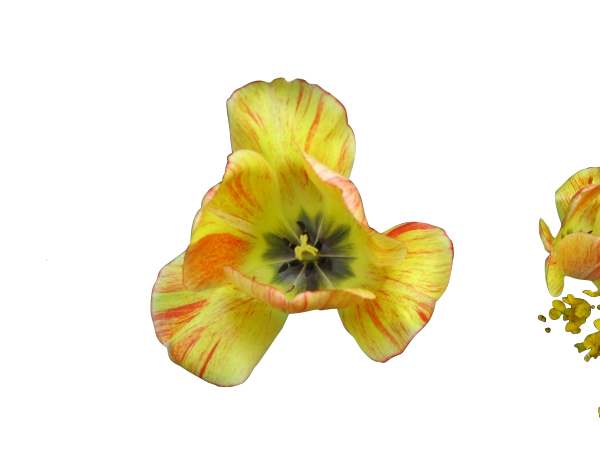}}
	\subfloat[Prior \label{fig:subfig:Yellowflowerprior}]{
		\includegraphics[width=1.8in, origin=br]{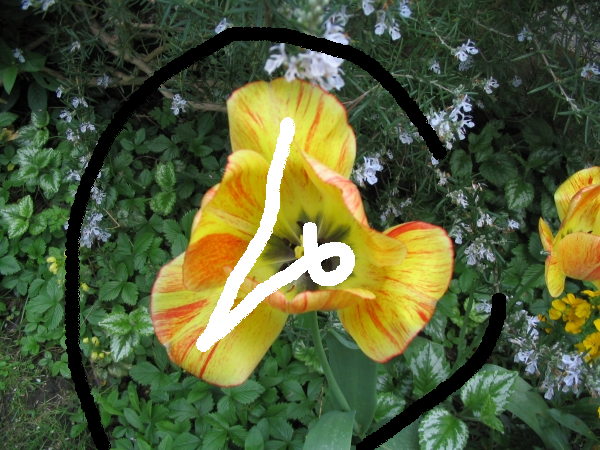}}
	\\
	\subfloat[pDCA final result (general method) \label{fig:subfig:Yellowflower_pDCA_914per_1294}]{
		\includegraphics[width=1.8in, origin=br]{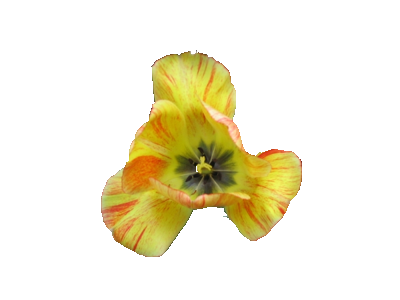}}
	\subfloat[svd1 final result \label{fig:subfig:Yellowflower_SVD_96per_135} ]{
		\includegraphics[width=1.8in, origin=br]{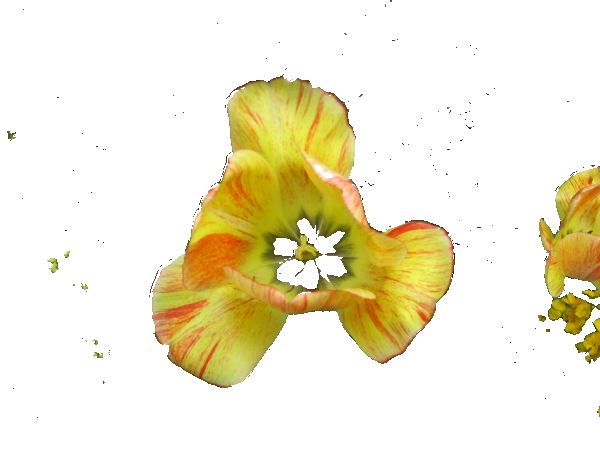}}
	\subfloat[pDCA final result (sparse method 1)\label{fig:yellowflower_pDCA_99per_593}]{
		\includegraphics[width=1.8in, origin=br]{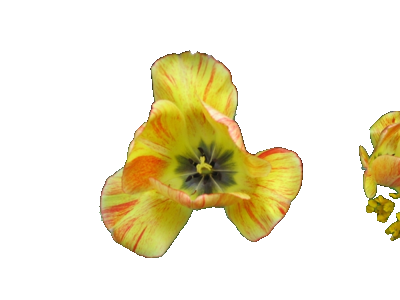}}
	\caption{Image {\bf{a}} is of size $400\times 300$. Image {\bf{b}} is the ground truth of segmented image {\bf{a}}. Image {\bf{c}} is the segmentation prior. Image {\bf{d}} is segmented by pDCA algorithm with a general search window in \ref{patch_and_window}b. Image {\bf{e}} is the result obtained by svd1 with $0.3\%$ sampling rate, and image {\bf{f}} is segmented images by pDCA algorithm with sparse method window in \ref{fig:search_window1}. }
	\label{fig:Segmentationbysparsemethod1} 
\end{figure*}

For the specially designed window \ref{fig:search_window1}, let's take the yellow flower \ref{fig:subfig:Originalyellowflower} for example. Due to the remote blocks of window \ref{fig:search_window1}, we can get more information in the far field. It can help us achieve better segmentation results for pDCA. It is shown in Figure \ref{fig:yellowflower_pDCA_99per_593} compared to the usual box search window segmented result as in Figure \ref{fig:subfig:Yellowflower_pDCA_914per_1294}. The small and disconnected yellow flower in the right part is missed in Figure \ref{fig:subfig:Yellowflower_pDCA_914per_1294}.



\begin{figure*}[!htbp]
	\centering
	\subfloat[Image: 400$\times$300 \label{fig:subfig:Originalstarfish}]{
		\includegraphics[width=1.8in, origin=br]{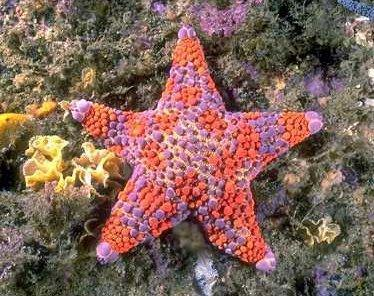}}
	\subfloat[Ground truth\label{fig:subfig:starfishgroundtruth}]{
		\includegraphics[width=1.8in, origin=br]{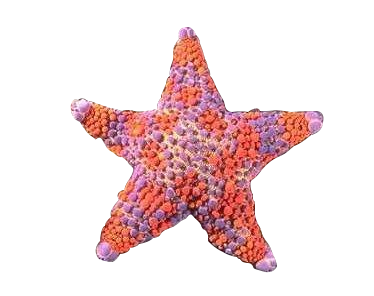}}
	\subfloat[Prior \label{fig:subfig:starfishprior}]{
		\includegraphics[width=1.8in, origin=br]{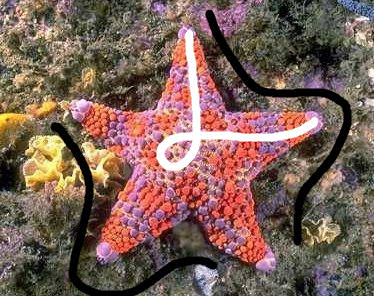}}
	\\
	\subfloat[pDCA final result (general
	method) \label{fig:subfig:starfish_pDCA_914per_1294}]{
		\includegraphics[width=1.8in, origin=br]{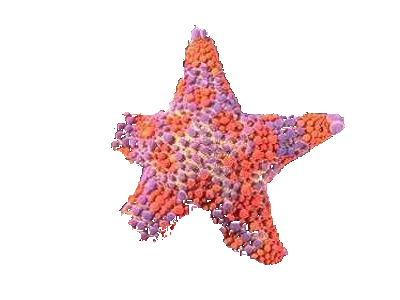}}
	\subfloat[svd1  result\label{fig:subfig:starfish_SVD_96per_135} ]{
		\includegraphics[width=1.8in, origin=br]{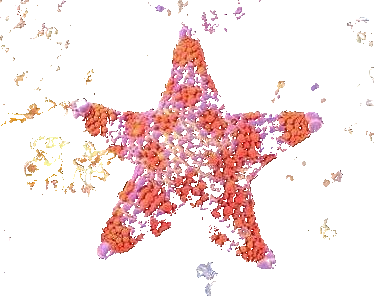}}
	\subfloat[pDCA final result (sparse
	method 2)\label{fig:subfig:starfish_pDCA_97per_404}]{
		\includegraphics[width=1.8in, origin=br]{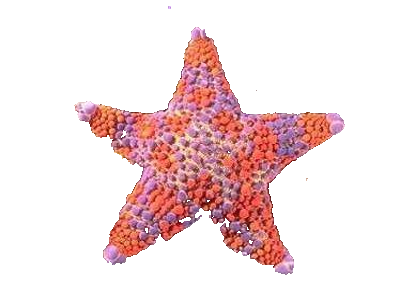}}
	\caption{Image {\bf{a}} contains $400\times 300$ pixels, image {\bf{b}} is the ideal segmentation result of image {\bf{a}}, image {\bf{c}} is the prior of this segmentation assignment, image {\bf{d}} is the result finished by pDCA algorithm with a general method which is shown in \ref{patch_and_window}b, image {\bf{e}} is the result finished by svd1 with $0.3\%$, and the image {\bf{f}} is the result finished by pDCA algorithm with sparse method 2 as showed in \ref{fig:search_window2}.}
	\label{fig:Segmentationbysparsemethod2} 
\end{figure*}

Now, let's turn to the starfish with sparse search window \ref{fig:search_window2} for example. There are many similar pixels in the background, which have quite similar colors to the body of starfish. It thus is difficult to separate the starfish from its background. As shown in Figures \ref{fig:subfig:starfish_pDCA_914per_1294} and \ref{fig:subfig:starfish_SVD_96per_135},  neither pDCA with general search window in \ref{patch_and_window}b nor svd1 with  $0.3\%$ sampling can segment the starfish well. However, with the specially designed search window \ref{fig:search_window2}, pDCA can obtain a promising segmentation in Figure \ref{fig:subfig:starfish_pDCA_97per_404}.  Compared with the general search window \ref{patch_and_window}b, the specially designed search window \ref{fig:search_window2} would have more information associated with the central pixel.  It can employ the information of some pixels which are far from the central pixel of the searching window. 

Now, let us turn to studying the efficiency of different preconditioners in numerics. 
\subsubsection{Comparison of different preconditioners}

In this section, we will compare the different preconditioners proposed in section \ref{subsec:preconditioner}. In our algorithm, there are three important factors for the performance of pDCA. They are the type of preconditioners, the type of Laplacian (normalized or unnormalized), and the type of step sizes (finite or infinite step sizes). Considering these factors, different results are shown under different conditions. We choose the image `stone' in figure \ref{fig:Segmentationbygeneralwindow_stone} for these comparisons. For the parameter in the proposed pDCA, the diffuse parameter $\epsilon$, the convex split parameter $c$, the fidelity parameter $\eta$, the search window size, and the number of precondition iterations are set to 100, 11, 100, $25\times 25$ and 4, respectively. According to the numerical experiments, the best choice for best performance is the normalized Laplacian operator with the infinity step size among a large number of combinations with the above factors.

As shown in Table \ref{tab:unnormal:normal:com}, for perturbed Jacobi preconditioner, the normalized Laplacian is better than the unnormalized graph Laplacian, since the normalized graph Laplacian is usually better conditioned. 
As shown in Table \ref{tab:normalized:lap:com}, for the normalized graph Laplacian, the perturbed Jacobi,  the damped Jacobi, and the generalized Richardson preconditioners are very competitive and efficient. Finally, as shown in Table \ref{tab:infinity:Jac} or \ref{tab:infinity:Ric}, for the proposed pDCA with damped Jacobi or generalized Richardson preconditioner, the step size infinity will bring the most efficient and stable algorithm. 
\begin{table}[h]
	\centering
	\tabcolsep = 0.5cm
	\begin{tabular}{|l|cc|c|}\hline
		method & \multicolumn{2}{c|}{perturbed Jacobi}&Stopping  \\
		&Normalized        &Unnormalized&Criterion\\\hline
		{time}            &\bf{54}                &106 &DICE=\\
		{iteration}       &\bf{282}               &724 &0.990\\\hline
	\end{tabular}
	\caption{Comparisons of different graph Laplacian operators using the perturbed Jacobi preconditioner.}
	\label{tab:unnormal:normal:com}
\end{table}
\begin{table}[h]
	\centering
	\tabcolsep = 0.5cm
	\begin{tabular}{|l|ccc|c|}\hline
		method   & \multicolumn{3}{c|}{ normalized graph Laplacian} &Stopping  \\
		&Perturbed&Damped&Richardson&Criterion\\
		\hline
		{time}            &57             &59            &61        &DICE=  \\
		{iteration}       &315            &330           &292        & 0.993\\\hline
	\end{tabular}
	\caption{Comparisons of different preconditioners using the normalized graph Laplacian. The `Perturbed', `Damped', and `Richardson' denote the perturbed Jacobi preconditioner, damped Jacobi preconditioner, and Richardson preconditioner correspondingly.} 
	\label{tab:normalized:lap:com}
\end{table}

\begin{table}[!htbp]
	\centering
	\tabcolsep = 0.45cm
	\begin{tabular}{|l|cccccc|c|}\hline
		method   & \multicolumn{6}{c|}{ normalized graph Laplacian, Perturbed Jacobi}&Stopping \\
		step size         &0.01    &0.05   &0.1    &1      &5      &$\infty$ &Criterion\\\hline
		{time}            &110     &68     &63     &59     &57     &57&DICE=\\
		{iteration}       &738     &390    &352    &318    &315    &315&0.993\\\hline
	\end{tabular}
	\caption{Comparisons of different step sizes using perturbed Jacobi preconditioners for the normalized graph Laplacian.}
	\label{tab:infinity:Jac}
\end{table}
\begin{table}[!htbp]
	\centering
	\tabcolsep = 0.45cm
	\begin{tabular}{|l|cccccc|c|}\hline
		method   & \multicolumn{6}{c|}{ normalized graph Laplacian, Richardson}&Stopping \\
		step size         &0.01    &0.05   &0.1    &1      &5      &$\infty$&Criterion\\\hline
		{time}            &112     &69     &64     &61     &60     &60&DICE=\\
		{iteration}       &721     &368    &328    &296    &292    &292&0.993\\\hline
	\end{tabular}
	\caption{Comparisons of different step sizes using Richardson preconditioners for the normalized graph Laplacian. }
	\label{tab:infinity:Ric}
\end{table}

\begin{figure*} [t!]
	\centering
	\includegraphics[width = 0.75\textwidth]{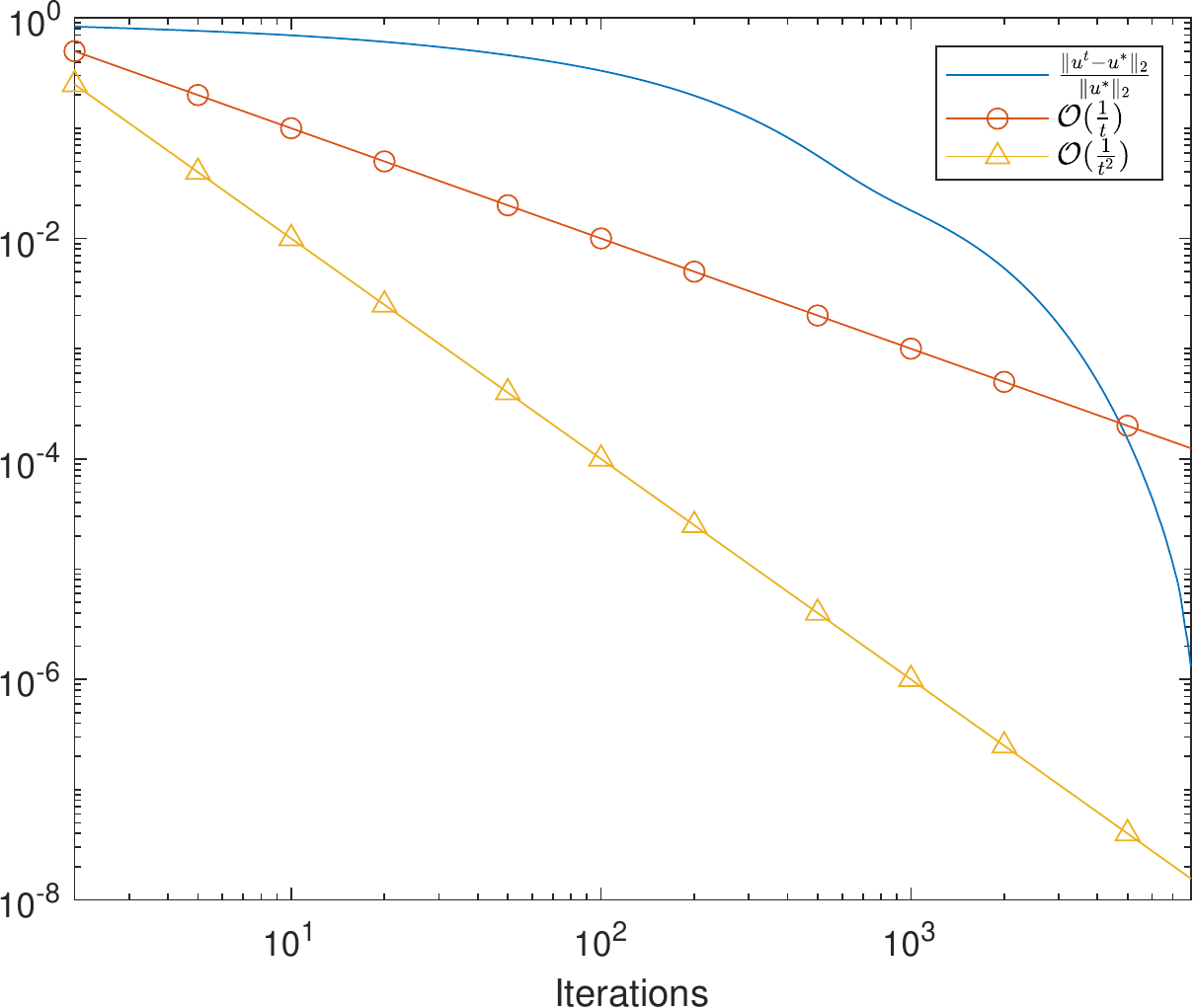}
	\caption{Numerical convergence rate of the iteration sequence $\{u^t\}$ with loglog plot where both the horizontal and vertical axes are changed to log scaling. The test image is `stone'. $u^*$ is obtained approximately from the iteration sequence $\{u^t\}$ after more than 10000 iterations, which is stable numerically.}
	\label{patch_and_window:rate} 
\end{figure*}

We now give some discussions on the local convergence rate of the proposed preconditioned DCA for image segmentation. As shown in Figure \ref{patch_and_window:rate}, the numerical convergence rate is better than the sublinear convergence rate as in Theorem \ref{thm:energy:decrease}.

Finally, let us turn to investigate the proposed pDCA for data clustering numerically. 
\subsection{Data Clustering}
This section will present some compact comparisons between the proposed pDCA algorithm and svd1 for data clustering  \cite{GMBFP}. The main tool, KNN \cite{MPL} for producing the weights for data clustering is implemented based on cluster tree data structure. It is not implemented in parallel although highly efficient. The proposed pDCA with KNN has fewer advantages than image segmentation. However, it is still very competitive compared to svd1.

We applied our algorithm to three data sets: MNIST (4 and 9), two moons, and a point set of concentric half circles. The MNIST data set is composed of 70,000 28 × 28 images of hand-written digits 0 through 9. We choose the combination of `4' and `9' which is more challenging to cluster than other combinations. To reduce the dimensions of data space, we transform the original data space to a smaller feature space with only $50$ dimensions by principal component analysis. Two moons data set was used by \cite{TMP} for spectral clustering with the graph $p$-Laplacian. It is constructed from the two half circles in $\mathbb{R}^2$ with radius one \cite{Tang2019}. High-dimension random noise is added to points in the two circles. 
The half-circle data set randomly chooses angles and two disjoint ranges of radius to generate a new data set. For the randomness of the algorithm, we calculate the mean of results in 10 experiments. The coefficient `Accuracy' indicates the proportion of points labeled as the correct classes. 

\begin{figure*} [t!]
	\centering
	\subfloat[MNIST49 by svd1\label{fig:mnist49bySVD}]{
		\includegraphics[width=1.8in,  origin=br]{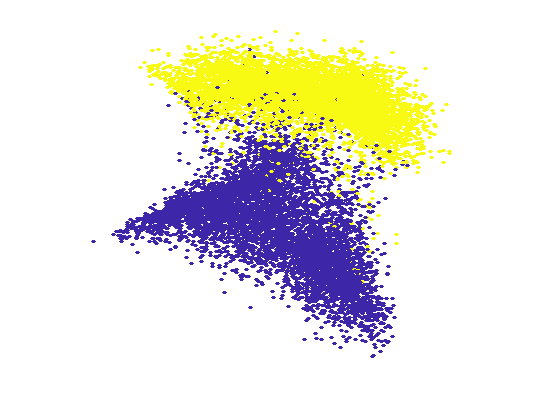}}\ \ 
	\subfloat[MNIST49 by pDCA\label{fig:mnist49bypDCA}]{
		\includegraphics[width=1.8in,  origin=br]{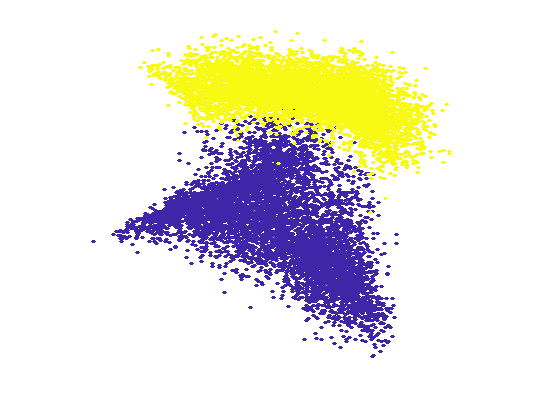}}\ \ 
	\caption{There are 13,782 data points in this combination, and we use the second and the third eigenvectors of the sparse graph Laplacian matrix to display the result of clustering in two dimensions. The parameters of the two algorithms are highly different. The parameters $\epsilon$, $c$, $\eta$, $dt$, and the number of eigenvectors are 1, 51, 50, 0.15, and 300 for svd1. The parameters $\epsilon$, $c$, $\eta$ are 100, 101,100 for pDCA.}
	\label{fig:clustering1:mnist} 
\end{figure*}

\begin{figure*} [t!]
	\centering
	\subfloat[two moons by svd1\label{fig:twomoonbySVD}]{
		\includegraphics[width=1.8in,  origin=br]{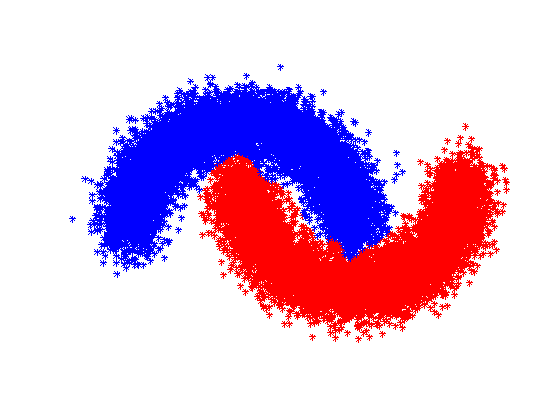}}\ \ 
	\subfloat[two moons by pDCA\label{fig:twomoonbypDCA}]{
		\includegraphics[width=1.8in,  origin=br]{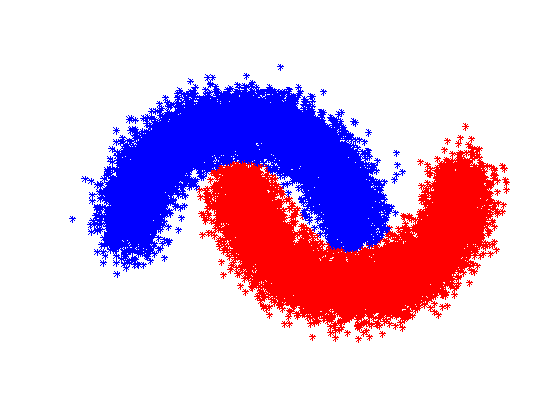}}\ \ 
	\caption{There are 20000 points in this data set. The parameters of the two algorithms are highly different. The parameters $\epsilon$, $c$, $\eta$, $dt$, and the number of eigenvectors are 1, 51, 50, 0.1, and 300 for svd1. The parameters $\epsilon$, $c$, $\eta$ are 100, 11,100 for pDCA.}
	\label{fig:dada:cluster:twomoons} 
\end{figure*}

\begin{figure*} [t!]
	\centering
	\subfloat[half-circle by svd1\label{fig:halfcirclebySVD}]{
		\includegraphics[width=1.8in,  origin=br]{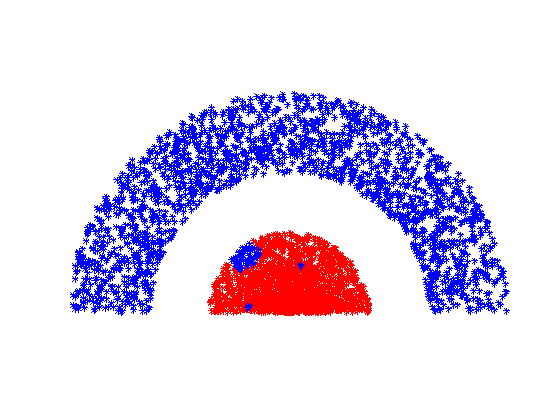}}\ \ 
	\subfloat[half-circle by pDCA\label{fig:halfcirclebypDCA}]{
		\includegraphics[width=1.8in,  origin=br]{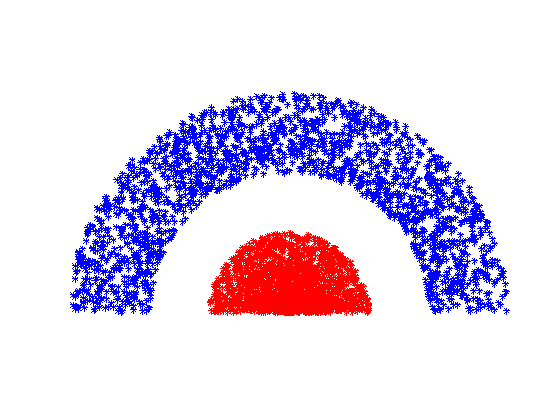}}\ \ 
	\caption{There are 4000 points in this data set. The parameters of the two algorithms are highly different. The parameters $\epsilon$, $c$, $\eta$, $dt$, and the number of eigenvectors are 1, 51, 50, 0.1, and 100 for svd1. The parameters $\epsilon$, $c$, $\eta$ are 100, 11,100 for pDCA.}
	\label{fig:cluster:half-circle} 
\end{figure*}

As shown in Table \ref{tab:data:clustering}, the proposed pDCA is very competitive compared with svd1 for the supervised case. We only present the results for pDCA for the unsupervised case which is very efficient since it is very challenging to find appropriate parameters for svd1. Figures \ref{fig:clustering1:mnist}, \ref{fig:dada:cluster:twomoons}, and \ref{fig:cluster:half-circle} show that the proposed pDCA is highly effective for data clustering with these data sets.


\begin{table}[h]
	\centering
	\begin{tabular}{|l|cc|cc|cc|cc|cc|}\hline
		&\multicolumn{6}{c|}{Supervised}&\multicolumn{4}{c|}{Unsupervised}\\\hline
		data set  & \multicolumn{2}{c|}{mnist49} & \multicolumn{2}{c|}{half-circle} & \multicolumn{2}{c|}{two moons} & \multicolumn{2}{c|}{half-circle} & \multicolumn{2}{c|}{two moons}\\
		&svd1        &pDCA    &svd1      &pDCA    &svd1     &pDCA&\multicolumn{2}{c|}{pDCA}&\multicolumn{2}{c|}{pDCA}\\ \hline
		time               &12.1      &11.2   &3.0   &2.3   &7.2    &3.1 &\multicolumn{2}{c|}{2.5}&\multicolumn{2}{c|}{3.1}\\
		iteration &100&100&1446&1&168 &4&\multicolumn{2}{c|}{1}&\multicolumn{2}{c|}{3}\\
		Accuracy          &0.9673    &\bf{0.9880} &0.9832  &\bf{1}  &0.9919 &\bf{0.9931} &\multicolumn{2}{c|}{{1}}&\multicolumn{2}{c|}{{0.9931}}\\ \hline
	\end{tabular}
	\caption{Comparisons between the proposed pDCA and  svd1 by iteration time and numbers. We employ the minst49, the half-circle, and the two moons data sets for comparisons for the supervised case. The parameters for svd1 and pDCA of the supervised case are in the corresponding captions of Figures \ref{fig:clustering1:mnist}, \ref{fig:dada:cluster:twomoons}, and \ref{fig:cluster:half-circle}.  We only present the performance of pDCA for the unsupervised case here. The parameters $\epsilon$, $c$, $\eta$ are 100, 11, and 100 for pDCA in the unsupervised case, and the parameters of KNN are $100$ for two moons and $10$ for half-circle. The performance of unsupervised clustering is surprisingly well especially for half-circle mainly because we use the second eigenvector of the normalized graph Laplacian matrix for initialization.}
	\label{tab:data:clustering}
\end{table}

\section{Conclusion}\label{sec:conclusion}
We mainly developed a  preconditioned DCA framework with parallel preconditioners for the graph Ginzburg-Landau model with applications for image segmentation and data clustering.  For the damped Jacobi or generalized Richardson preconditioned iteration, the main computation is the matrix-vector multiplication. The NFFT (Nonequispaced fast Fourier transform) can employ the structure of Gaussian kernel \cite{BSV} and may help accelerate the matrix-vector multiplication on GPU. For date clustering, the parallel implementation of KNN will bring out great benefits for data clustering in the proposed framework. Besides, more general nonlocal energy functional in \cite{BBTB} is also very interesting within the proposed pDCA framework.




\noindent
{\small
	\textbf{Acknowledgements}
	Xinhua Shen and Hongpeng Sun acknowledge the support of the National Natural Science Foundation of China under grant No. \,12271521 and Beijing Natural Science Foundation No. Z210001. The work of Xuecheng Tai was supported by RG(R)-RC/17-18/02-MATH, HKBU 12300819, NSF/RGC Grant N-HKBU214-19  and RC-FNRA-IG/19-20/SCI/01.
}


\end{document}